\documentclass[11pt]{article}
\usepackage{amsmath, amssymb, amsfonts, amsthm}

\newtheorem{theorem}{Theorem}

\newtheorem{proposition}{Proposition}

\newtheorem{lemma}{Lemma}

\theoremstyle{remark}
\newtheorem{remark}{Remark}
\begin{document}

\title{\textbf{Couplings and Strong Approximations to Time Dependent
Empirical Processes Based on I.I.D. Fractional Brownian Motions}}
\author{P\'eter Kevei \\
MTA-SZTE Analysis and Stochastics Research Group, Bolyai Institute \\
Aradi v\'{e}rtan\'{u}k tere 1, 6720 Szeged, Hungary, and \\
Center for Mathematical
Sciences, Technische Universit\"at M\"unchen \\
Boltzmannstra{\ss }e 3, 85748
Garching, Germany \\
e-mail: \texttt{kevei@math.u-szeged.hu} 
\and David M. Mason \\
Department of Applied Economics and Statistics, University of Delaware \\
213 Townsend Hall, Newark, DE 19716, USA \\
e-mail: \texttt{davidm@udel.edu}}
\date{}
\maketitle

\begin{abstract}
\noindent We define a time dependent empirical process based on $n$
i.i.d.~fractional Brownian motions and establish Gaussian couplings and
strong approximations to it by Gaussian processes. They lead to functional
laws of the iterated logarithm for this process.

 \smallskip

\end{abstract}

\section{Introduction}

\label{s1}

The aim in this paper is to derive Gaussian couplings and strong
approximations to time dependent empirical processes based on $n$
independent sample continuous fractional Brownian motions, as defined in
Subsection \ref{s2A}. Our couplings yield surprisingly close \textit{almost
sure} approximations of our empirical processes by Gaussian processes
defined on sequences of intervals for which weak convergence cannot hold in
the limit. As an example of what our strong approximations can do, we show
that functional laws of the iterated logarithm [FLIL] for these empirical
processes can be derived from those that are known for Gaussian
processes.\smallskip

Our investigations may be thought of as a continuation of those of Kuelbs,
Kurtz and Zinn \cite{KKZ}, who proved central limit theorems for time
dependent empirical processes based on $n$ independent copies of a wide
variety of random processes. These include certain self-similar processes of
which fractional Brownian motion is a special case. Our results reveal the
kind of strong limit theorems that are possible when one turns to the
detailed analysis of time dependent empirical processes based on processes
which have a fine local random structure, such as fractional Brownian
motion.\smallskip

Kuelbs and Zinn \cite{KZ1,KZ2} have obtained central limit theorems for a
time dependent quantile process based on $n$ independent copies of a wide
variety of random processes. In the process they generalized a result of
Swanson \cite{Swan07}, who used classical weak convergence theory to prove
that an appropriately scaled median of $n$ independent Brownian motions
converges weakly to a mean zero Gaussian process. In a sequel to this paper
we use the results in the present work to derive strong approximations and
FLILs for quantile processes or inverses of these time dependent empirical
processes based on $n$ i.i.d.~sample continuous fractional Brownian motions.
For details see Kevei and Mason \cite{KM}. \smallskip

To motivate our work, we point out some implications of a coupling and a
strong approximation due to Koml\'{o}s, Major and Tusn\'{a}dy (KMT) \cite%
{KMT}. Let $X_{1},X_{2},\dots,$ be i.i.d.~$F$. For each $n\geq1$, let 
\begin{equation*}
F_{n}\left( x\right) =n^{-1}\sum_{j=1}^{n}1\left\{ X_{j}\leq x\right\}, \
x\in\mathbb{R},
\end{equation*}
denote the empirical distribution function based on $X_{1},\dots,X_{n}$, and
define the empirical process 
\begin{equation*}
v_{n}\left( x\right) :=\sqrt{n}\left\{ F_{n}\left( x\right) -F\left(
x\right) \right\} ,\ x\in\mathbb{R}.
\end{equation*}
Using the coupling result given in Theorem 3 of KMT \cite{KMT} one can
construct a probability space on which sit an i.i.d.~$F$ sequence $X_{1}$, $%
X_{2}$, $\ldots$, and a sequence of Brownian bridges $B_{1},B_{2},\ldots$,
on $[0,1]$ such that 
\begin{equation}
\left\Vert v_{n} -B_{n}\left( F\right) \right\Vert _{\mathbb{R}}=O\left( 
\frac{\log n}{\sqrt{n}}\right) ,\text{ a.s.,}  \label{cup}
\end{equation}
where for a real-valued function $\Upsilon$ defined on a set $S$ we use the
notation 
\begin{equation}
\left\Vert \Upsilon\right\Vert _{S}=\sup_{s\in S}\left\vert \Upsilon\left(
s\right) \right\vert .  \label{sup}
\end{equation}
The rate $\log n/\sqrt{n}$ in (\ref{cup}) is optimal.\smallskip

Further, by the strong approximation result stated in Theorem 4 of KMT \cite%
{KMT} one has on the same probability space an i.i.d.~$F$ sequence $X_{1}$, $%
X_{2}$, $\ldots,$ and a sequence of \textit{independent} Brownian bridges $%
B_{1},B_{2},\ldots$, on $[0,1]$ such that 
\begin{equation}
\left\Vert v_{n}-\frac{\sum_{j=1}^{n}B_{j}\left( F\right) }{\sqrt{n}}%
\right\Vert _{\mathbb{R}}=O\left( \frac{\left( \log n\right) ^{2}}{\sqrt{n}}%
\right) ,\text{ a.s.}  \label{kmt}
\end{equation}
It is known that the $n^{-1/2}$ part of the rate in (\ref{kmt}) is optimal,
but not the $\left( \log n\right) ^{2}$. It has long been conjectured that
the $\left( \log n\right) ^{2}$ in (\ref{kmt}) can be replaced by $\log n$.
This is one of the rare cases where any such optimality is known in the rate
of strong approximation to an empirical process.

Our goal is to develop analogs of (\ref{cup}) and (\ref{kmt}) for the time
dependent empirical processes based on independent copies of sample
continuous fractional Brownian motion. These are described in the next
section. The rates of coupling and strong approximation that we obtain are
unlikely to be anywhere near optimal in the sense just described, however
they will be seen to be sufficient to derive from them FLILs for our time
dependent empirical processes. We find it noteworthy that useful couplings
and strong approximations can be obtained for the kind of complexly formed
empirical processes that we consider. Our main results are detailed in
Section \ref{s2} and they are proved in Section \ref{s4}. We gather together
some needed facts in the Appendix. To prove our main results we use the
methodology outlined in Berthet and Mason \cite{BerthetMason}.

\section{Coupling and strong approximation to a time dependent empirical
process}
\label{s2}

\subsection{A time dependent empirical process}
\label{s2A}

\noindent Let $\big\{ B^{(H)}\big\} \cup \big\{ B_{j}^{(H)}\big\}
_{j\geq 1}$ be a sequence of i.i.d.~sample continuous fractional Brownian
motions with Hurst index $0<H<1$ defined on $[0,\infty )$. Note that $%
B^{(H)} $ is a continuous mean zero Gaussian process on $[0,\infty )$ with
covariance function defined for any $s,t\in \lbrack 0,\infty )$%
\begin{equation*}
E\left( B^{(H)}\left( s\right) B^{(H)}\left( t\right) \right) =\frac{1}{2}%
\left( \left\vert s\right\vert ^{2H}+\left\vert t\right\vert
^{2H}-\left\vert s-t\right\vert ^{2H}\right).
\end{equation*}%
By the L\'{e}vy modulus of continuity theorem for sample continuous
fractional Brownian motion $B^{(H)}$ with Hurst index $0<H<1$ (see (\ref{W})
below), we have for any $0<T<\infty $, w.p.~$1$, 
\begin{equation}
\sup_{0\leq s\leq t\leq T}\frac{\left\vert B^{(H)}\left( t\right)
-B^{(H)}\left( s\right) \right\vert }{f_{H}(t-s)}=:L<\infty ,  \label{MC}
\end{equation}%
where for $u\geq 0$ 
\begin{equation}
f_{H}(u)=u^{H}\sqrt{1\vee \log u^{-1}}  \label{FH}
\end{equation}%
and $a\vee b=\max \{a,b\}$. We shall take versions of $\big\{B^{(H)}\big\} \cup \big\{ 
B_{j}^{(H)}\big\}_{j\geq 1}$ such that (\ref{MC}) holds for all of their trajectories. 
\smallskip

For any $t\in \left[ 0,\infty \right) $ and $x\in \mathbb{R}$ let $F\left(
t,x\right) =P\left\{ B^{(H)}\left( t\right) \leq x\right\} .$ Note that 
\begin{equation}
F\left( t,x\right) =\Phi \left( x/t^{H}\right) ,  \label{Ft}
\end{equation}%
where $\Phi \left( x\right) =P\left\{ Z\leq x\right\} ,$ with $Z$ being a
standard normal random variable. For any $n\geq 1$ define the time dependent 
\textit{empirical distribution function} 
\begin{equation*}
F_{n}\left( t,x\right) =n^{-1}\sum_{j=1}^{n}1\left\{ B_{j}^{\left( H\right)
}\left( t\right) \leq x\right\} .
\end{equation*}
Applying Theorem 5 in \cite{KKZ} (also see their Remark 8) one can
show for any choice of $0<\gamma \leq 1<T<\infty $ that the time dependent 
\textit{empirical process} indexed by $\left( t,x\right) \in \mathcal{T}%
\left( \gamma \right) $, 
\begin{equation*}
v_{n}\left( t,x\right) =\sqrt{n}\left\{ F_{n}\left( t,x\right) -F\left(
t,x\right) \right\} ,
\end{equation*}%
where%
\begin{equation*}
\mathcal{T}\left( \gamma \right) :=\left[ \gamma ,T\right] \times \mathbb{R},
\end{equation*}%
converges weakly to a uniformly continuous centered Gaussian process $%
G\left( t,x\right) $ indexed by $\left( t,x\right) \in \mathcal{T}\left(
\gamma \right) $, whose trajectories are bounded, having covariance function%
\begin{equation}
\begin{split}
& E\left( G\left( s,x\right) G\left( t,y\right) \right) \\
& =P\left\{ B^{\left( H\right) }\left( s\right) \leq x,B^{(H)}\left(
t\right) \leq y\right\} -P\left\{ B^{(H)}\left( s\right) \leq x\right\}
P\left\{ B^{(H)}\left( t\right) \leq y\right\} .
\end{split}
\label{EG}
\end{equation}%
Keeping in mind that $\mathcal{T}\left( \gamma \right) $ is equipped with
the semimetric

\begin{equation}
\rho \left( (s,x\right) ,(t,y))=\sqrt{E\left( G\left( s,x\right) -G\left(
t,y\right) \right) ^{2}},  \label{metric}
\end{equation}%
we see by weak convergence that $\mathcal{T}\left( \gamma \right)$ is
totally bounded and thus separable in the topology induced by this
semimetric $\rho$. Moreover its completion $\mathcal{T}^{c}\left( \gamma
\right) $ in this topology is compact. Since $G$ is bounded and uniformly
continuous on $\mathcal{T}\left( \gamma \right) $ it can be extended
uniquely to be bounded and uniformly continuous on $\mathcal{T}^{c}\left(
\gamma \right) $.

\begin{remark}
\label{rem:c} To see how this is done, notice that for each $t\in \left[
\gamma ,T\right]$, both $\left\{ \left( t,-m\right) \right\} _{m\geq 1}$ and 
$\left\{ \left( t,m\right) \right\} _{m\geq 1}$ are Cauchy sequences in $%
\mathcal{T}\left( \gamma \right) $ with respect to the semimetric $\rho $. \
Also by the boundedness and uniform continuity of $G$ on $\mathcal{T}\left(
\gamma \right) $, the sequences $\left\{ G\left( t,-m\right) \right\}
_{m\geq 1}$ and $\left\{ G\left( t,m\right) \right\} _{m\geq 1}$ are also
bounded Cauchy sequences in $\mathbb{R}$. Furthermore, both $EG^{2}\left(
t,-m\right) \rightarrow 0$\ and $EG^{2}\left( t,m\right) \rightarrow 0$, as $%
m\rightarrow \infty $. Thus we can unambiguously define $\left( t,-\infty
\right) $ as the limit of the sequence $\left( t,-m\right) $ as $%
m\rightarrow \infty $ and $G\left( t,-\infty \right) =0$, w.p.~$1$, and $%
\left( t,\infty \right) $ as the limit of the sequence $\left( t,m\right) $
as $m\rightarrow \infty $ and $G\left( t,\infty \right) =0$, w.p.~$1$. We
see that for any $t\in \left[ \gamma ,T\right] $ and $\left( s,y\right) \in 
\mathcal{T}\left( \gamma \right)$, 
\begin{equation*}
\rho \left( \left( t,\pm \infty \right) ,\left( s,y\right) \right) =\sqrt{%
E\left( G(\left( t,\pm \infty \right) )-G\left( s,y\right) \right) ^{2}}=%
\sqrt{EG^{2}\left( s,y\right) }
\end{equation*}%
and for $s,t\in \left[ \gamma ,T\right] $ 
\begin{equation*}
\rho \left( \left( t,\pm \infty \right) ,\left( s,\pm \infty \right) \right)
=\sqrt{E\left( G(\left( t,\pm \infty \right) )-G\left( s,\pm \infty \right)
\right) ^{2}}=0.
\end{equation*}%
With these definitions $\rho $ becomes a semimetric on 
$\left[ \gamma ,T\right] \times \left( \mathbb{R\cup }
\left\{ -\infty ,\infty \right\} \right)$. Next consider 
$\left[ \gamma,T\right] \times \left\{ -\infty,\infty \right\}$
as an equivalence class, i.e.~$\left( t,\pm \infty\right) \sim
\left( s,\pm \infty \right)$, whenever
$\rho \left( \left(t,\pm \infty \right) ,
\left( s,\pm \infty \right) \right) =0$, which always
happens, and denote it by $\omega $ and with some abuse of the previous
notation write $G\left( \omega \right) =0$, $\rho \left( \omega ,\omega
\right) =0$ and for any $\left( s,y\right) \in \mathcal{T}\left( \gamma
\right) $, $\rho \left( \omega ,\left( s,y\right) \right) =\sqrt{%
EG^{2}\left( s,y\right) },$ and let $\rho $ remain as it was previously
defined on $\mathcal{T}\left( \gamma \right) \times \mathcal{T}\left( \gamma
\right) $. We define the completion of $\mathcal{T}^{c}\left( \gamma \right)
= \left( \left[ \gamma ,T\right] \times \mathbb{R} \right) \cup \left\{
\omega \right\}$, which is readily shown to be a complete metric space with
semimetric $\rho$.
\end{remark}

Therefore we can consider $G$ as a Gaussian process taking values in the
separable Banach space consisting of the continuous functions in the
sup-norm on the compact metric space $\mathcal{T}^{c}\left( \gamma \right)$.
For later use we point out that by Proposition 1 on page 26 of Lifshits \cite%
{Lifshits} we can assume that the Gaussian process $G\left( t,x\right) $ is
separable.\smallskip

For future reference we record here that for some finite positive constant $%
M\left( \gamma ,T,H\right) $ for all $n\geq 1$ 
\begin{equation}
E\left\Vert v_{n}\right\Vert _{\mathcal{T}\left( \gamma \right) }\leq
M\left( \gamma ,T,H\right) .  \label{M}
\end{equation}%
Assertion (\ref{M}) follows from an application of the Hoffmann--J\o rgensen
inequality, cf.~Ledoux and Talagrand \cite{LedeouxTalagrand3}, page 156. For
the argument see, for instance, Lemma 3.1 of Einmahl and Mason \cite{EM97}.
\smallskip

We restrict ourselves to positive $\gamma $, since in Section 8.1 of
\cite{KKZ} it is pointed out that the empirical process 
$v_{n}\left( t,x\right)$
indexed by $\mathcal{T}\left( 0\right) :=\left[ 0,T\right] \times \mathbb{R}$
does not converge weakly to a uniformly continuous centered Gaussian process
indexed by $\left( t,x\right) \in \mathcal{T}\left( 0\right) $, whose
trajectories are bounded. More generally in the sequel, $G\left( t,x\right)$
denotes a centered Gaussian process on $\mathcal{T}\left( 0\right) $ with
covariance (\ref{EG}) that is uniformly continuous on
$\mathcal{T}\left(\gamma \right)$ with bounded trajectories for any
$0<\gamma \leq 1<T<\infty$.

We shall also be using the following empirical process indexed by function
notation. Let $X,X_{1},X_{2},\dots $, be i.i.d.~random variables from a
probability space $\left( \Omega ,\mathcal{A},P\right) $ to a measurable
space $\left( S,\mathcal{S}\right) $. Consider an empirical process indexed
by a class $\mathcal{G}$ of bounded measurable real valued functions on $%
\left( S,\mathcal{S}\right) $ defined by 
\begin{equation*}
\alpha _{n}\left( \varphi \right) :=\sqrt{n}(P_{n}-P)\varphi =\frac{%
\sum_{i=1}^{n}\varphi \left( X_{i}\right) -nE\varphi \left( X\right) }{\sqrt{%
n}}\text{, }\varphi \in \mathcal{G},
\end{equation*}%
where 
\begin{equation*}
P_{n}\left( \varphi \right) =n^{-1}\sum_{i=1}^{n}\varphi \left( X_{i}\right) 
\text{ and }P\left( \varphi \right) =E\varphi \left( X\right) \text{.}
\end{equation*}%
Keeping this notation in mind, let $\mathcal{C}\left[ 0,T\right] $ be the
class of continuous functions $g$ on $\left[ 0,T\right] $ endowed with the
topology of uniform convergence and where $\mathcal{B}\left[ 0,T\right] $
denotes the Borel subsets of $\mathcal{C}\left[ 0,T\right] $. Define this
subclass of $\mathcal{C}\left[ 0,T\right]$ 
\begin{equation}
\mathcal{C}_{\infty }:=\left\{ g:\ \sup \left\{ \frac{\left\vert g\left(
s\right) -g\left( t\right) \right\vert }{f_{H}(\left\vert s-t\right\vert )}%
,\ 0\leq s,t\leq T\right\} <\infty \right\} .  \label{CI}
\end{equation}%
Further, let $\mathcal{F}_{\left( \gamma ,T\right) }$ be the class of
functions of $g\in $ $\mathcal{C}\left[ 0,T\right] \rightarrow \mathbb{R}$,
indexed by $\left( t,x\right) \in \mathcal{T}\left( \gamma \right) ,$ of the
form 
\begin{equation*}
h_{t,x}\left( g\right) =1\left\{ g\left( t\right) \leq x,g\in \mathcal{C}%
_{\infty }\right\} .
\end{equation*}%
Here we permit $\gamma =0$. Since by (\ref{MC}) we can assume that each $%
B^{(H)},$ $B_{j}^{(H)}$, $j\geq 1$, is in $\mathcal{C}_{\infty }$, we see
that for any $h_{t,x}\in \mathcal{F}_{\left( \gamma ,T\right) }$, 
\begin{equation*}
\alpha _{n}\left( h_{t,x}\right) =\frac{1}{\sqrt{n}}\sum_{i=1}^{n}\left(
1\left\{ B_{i}^{(H)}\left( t\right) \leq x\right\} -P\left\{ B^{(H)}\left(
t\right) \leq x\right\} \right) =v_{n}\left( t,x\right) .
\end{equation*}%
We shall be using the notation $\alpha _{n}\left( h_{t,x}\right) $ and $%
v_{n}\left( t,x\right) $ interchangeably. \smallskip

Let $\mathbb{G}_{\left( \gamma ,T\right) }$ denote the mean zero Gaussian
process indexed by $\mathcal{F}_{\left( \gamma ,T\right) }$, having
covariance function defined for $h_{s,x},h_{t,y}\in $ $\mathcal{F}_{\left(
\gamma ,T\right) }$ 
\begin{equation*}
E\left( \mathbb{G}_{\left( \gamma ,T\right) }\left( h_{s,x}\right) \mathbb{G}%
_{\left( \gamma ,T\right) }\left( h_{t,y}\right) \right) =P\left\{ B^{\left(
H\right) }\left( s\right) \leq x,B^{(H)}\left( t\right) \leq y,B^{\left(
H\right) }\in \mathcal{C}_{\infty }\right\}
\end{equation*}%
\begin{equation*}
-P\left\{ B^{(H)}\left( s\right) \leq x,B^{(H)}\in \mathcal{C}_{\infty
}\right\} P\left\{ B^{(H)}\left( t\right) \leq y,B^{(H)}\in \mathcal{C}%
_{\infty }\right\} ,
\end{equation*}%
which since $P\left\{ B^{(H)}\in \mathcal{C}_{\infty }\right\} =1$,%
\begin{equation*}
=E\left( G\left( s,x\right) G\left( t,y\right) \right) ,
\end{equation*}%
i.e.~$\mathbb{G}_{\left( \gamma ,T\right) }\left( h_{t,x}\right) $ defines a
probabilistically equivalent version of the Gaussian process $G\left(
t,x\right) $ for $\left( t,x\right) \in \mathcal{T}\left( \gamma \right) $.
We shall say that a process $\widetilde{\mathcal{Y}}$ is a \textit{%
probabilistically equivalent version} of $\mathcal{Y}$ if $\widetilde{%
\mathcal{Y}}\overset{\mathrm{D}}{=}\mathcal{Y}$. \smallskip

Notice that in this notation 
\begin{equation*}
\begin{split}
\rho \left( (s,x\right) ,(t,y)) & =\sqrt{E\left( \mathbb{G}_{\left( \gamma
,T\right) }\left( h_{s,x}\right) -\mathbb{G}_{\left( \gamma ,T\right)
}\left( h_{t,y}\right) \right) ^{2}} \\
& =\sqrt{\mathop{Var} \left( h_{s,x}\left( B^{\left( H\right) }\right)
-h_{t,y}\left( B^{\left( H\right) }\right) \right) } \\
& \leq \sqrt{E\left( h_{s,x}\left( B^{\left( H\right) }\right)
-h_{t,y}\left( B^{\left( H\right) }\right) \right) ^{2}}=:d_{P}\left(
h_{s,x},h_{t,y}\right) .
\end{split}%
\end{equation*}%
More generally, for suitable functions $f$ and $g$ we shall write%
\begin{equation}
d_{P}\left( f,g\right) =\sqrt{E\left( f\left( B^{\left( H\right) }\right)
-g\left( B^{\left( H\right) }\right) \right) ^{2}}.  \label{DP}
\end{equation}

The proofs of a number our results rely on a lemma of Berkes and Philipp 
\cite{BerkesPhilipp}, which for the convenience of the reader we state
here.\smallskip

\noindent\textbf{Lemma A1 of Berkes and Philipp (1979)}
\textit{Let $S_{i},i=1,2,3$ be Polish spaces. Let
$\mathbf{F}$ be a distribution on $S_{1}\times S_{2}$ and 
$\mathbf{G}$ be a distribution on $S_{2}\times S_{3}$ such that the 
second marginal of $\mathbf{F}$ is equal to the first marginal of 
$\mathbf{G}$. Then there exists a probability space and a random vector 
$\left( Z_{1}, Z_{2}, Z_{3} \right)$ defined on it taking its values in
$S_{1}\times S_{2}\times S_{3}$ such that $\left( Z_{1},Z_{2}\right)$ 
has distribution $\mathbf{F}$ and $\left( Z_{2},Z_{3} \right)$
has distribution $\mathbf{G}$.}

\subsection{Our main coupling and strong approximation results for $\protect%
\alpha _{n}$}
\label{main}

In the results that follow 
\begin{equation}
\nu _{0}=2+\frac{2}{H}\quad \text{and }\ H_{0}=1+H.  \label{nuk}
\end{equation}%
We have the following Gaussian coupling to the empirical process $\alpha
_{n} $ indexed by $\mathcal{F}_{(\gamma _{n},T)}$.

\begin{proposition}
\label{prop:1} As long as $0<\gamma _{n}\leq 1$ satisfies for some $0\leq
\eta <\frac{1}{5H_{0}}$, 
\begin{equation}
\infty >-\frac{\log \gamma _{n}}{\log n}\rightarrow \eta ,\ \text{ as }%
n\rightarrow \infty ,  \label{eta}
\end{equation}%
for every $\lambda >1$ there exists a $\rho (\lambda )>0$ such that for each
integer $n$ large enough one can construct on the same probability space
random vectors $B_{1}^{(H)}$, $\ldots $, $B_{n}^{(H)}$ i.i.d.~$B^{(H)}$ and
a probabilistically equivalent version $\widetilde{\mathbb{G}}_{(\gamma
_{n},T)}^{(n)}$ of $\mathbb{G}_{(\gamma _{n},T)}$ such that, 
\begin{equation}
P\left\{ \left\Vert \alpha _{n}-\widetilde{\mathbb{G}}_{(\gamma
_{n},T)}^{\left( n\right) }\right\Vert _{\mathcal{F}_{(\gamma _{n},T)}}>\rho
(\lambda )\left( \log n\right) ^{\tau _{2}}\left( n^{-1/2}\gamma
_{n}^{-5H_{0}/2}\right) ^{2/(2+5\nu _{0})}\right\} \leq n^{-\lambda },
\label{p1}
\end{equation}%
where $\tau _{2}=(19H+25)/(24H+20)$ and $\nu _{0}$ is defined in (\ref{nuk}).
Moreover, in particular, when $\gamma _{n}=n^{-\eta }$, with
$0\leq \eta < \frac{1}{5H_{0}}$, 
\begin{equation*}
P\left\{ \left\Vert \alpha _{n}-\mathit{\ }\widetilde{\mathbb{G}}_{\left(
\gamma _{n},T\right) }^{(n)}\right\Vert _{\mathcal{F}_{(\gamma
_{_{n}},T)}}>\rho \left( \lambda \right) n^{-\tau _{1}}\left( \log n\right)
^{\tau _{2}}\right\} \leq n^{-\lambda },
\end{equation*}%
where $\tau _{1}=\tau _{1}(\eta )=\left( 1-5H_{0}\eta \right) /\left( 2+5\nu
_{0}\right) >0$.
\end{proposition}

\begin{remark} \label{rem:1}
Notice that Proposition \ref{prop:1} yields the coupling rate 
\begin{equation}
\left\Vert \alpha _{n}-\mathit{\ }\widetilde{\mathbb{G}}_{\left( \gamma
_{n},T\right) }^{(n)}\right\Vert _{\mathcal{F}_{(\gamma
_{n},T)}}=O_{P}\left( \left( \log n\right) ^{\tau _{2}}\left( n^{-1/2}\gamma
_{n}^{-5H_{0}/2}\right) ^{2/(2+5\nu _{0})}\right) .  \label{ppp}
\end{equation}%
In particular, for any $0<H<1$ and $0<\eta <1/\left( 5H_{0}\right)$ the
convergence (\ref{ppp}) holds with $\gamma _{n}=n^{-\eta }$, since such
$\gamma_{n}$ satisfy (\ref{eta}). The convergence (\ref{ppp}) is surprising
in light of the results in Section 8.1 in \cite{KKZ}, where it is pointed
out that the empirical process $v_{n}\left( t,x\right) $ indexed by
$\left[0,T\right] \times \mathbb{R}$ does not converge weakly to a uniformly
continuous centered Gaussian process indexed by
$\left( t,x\right) \in \left[0,T\right] \times \mathbb{R}$
whose trajectories are bounded. Observe,
however, by Theorem 5 in \cite{KKZ} for each $n\geq 1$ there is a version of
Gaussian process $G_{n}\left( t,x\right) =
\widetilde{\mathbb{G}}_{\left( \gamma _{n},T\right) }\left( h_{t,x}\right)$, 
which is uniformly continuous
on $\left[ \gamma _{n},T\right] \times \mathbb{R}$ with bounded
trajectories. We shall see that a coupling result following from a special
case of Theorem 1.1 of Zaitsev \cite{87a} is crucial to establish (\ref{p1})
on intervals $\left[ \gamma _{n},T\right] $ such that $\gamma _{n}$ goes to
zero at the rate (\ref{eta}).
\end{remark}

For any $\kappa>0$ let 
\begin{equation}
\mathcal{G}\left( \kappa\right) =\left\{ t^{\kappa}h_{t,x}:\left( t,x\right)
\in\left[ 0,T\right] \times\mathbb{R}\right\} .  \label{Gk}
\end{equation}
For $g\in$ $\mathcal{G}\left( \kappa\right)$, with some abuse of notation,
we shall write%
\begin{equation}
\mathbb{G}_{\left( 0,T\right) }\left( g\right) =t^{\kappa}\mathbb{G}_{\left(
0,T\right) }\left( h_{t,x}\right) .  \label{not}
\end{equation}
Also, in analogy with (\ref{sup}), we set 
\begin{equation*}
\left\Vert \alpha_{n}-\mathbb{G}_{\left( 0,T\right) }^{\left( n\right)
}\right\Vert _{\mathcal{G}\left( \kappa\right) }\!:=\sup\left\{ \left\vert
t^{\kappa}\alpha_{n}\left( h_{t,x}\right) -t^{\kappa}\mathbb{G}_{\left(
0,T\right) }^{(n)}\left( h_{t,x}\right) \right\vert :\left( t,x\right) \in%
\left[ 0,T\right] \times\mathbb{R}\right\}.
\end{equation*}
We get the following Gaussian coupling to the empirical process $\alpha_{n}$
indexed by $\mathcal{G}\left( \kappa\right) $.

\begin{proposition} \label{prop:2}
For any $0<\kappa <\infty $ and every 
$\lambda >1$ there
exists a $\rho ^{\prime }\left( \lambda \right) >0$ such that for each
integer $n$ large enough one can construct on the same probability space
random vectors $B_{1}^{(H)},\ldots ,B_{n}^{(H)}$ i.i.d.~$B^{\left( H\right)
} $ and a probabilistically equivalent version 
$\widetilde{\mathbb{G}}_{\left( 0,T\right) }^{(n)}$ of
$\mathbb{G}_{\left( 0,T\right) }$ such that, 
\begin{equation}
P\left\{ \left\Vert \alpha _{n}-\widetilde{\mathbb{G}}_{\left( 0,T\right)
}^{(n)}\right\Vert _{\mathcal{G}\left( \kappa \right) }>\rho ^{\prime
}\left( \lambda \right) n^{-\tau _{1}^{\prime }}\left( \log n\right) ^{\tau
_{2}}\right\} \leq n^{-\lambda },  \label{del}
\end{equation}
where $\tau _{2}$ is as in Proposition \ref{prop:1} and
$\tau _{1}^{\prime}=\tau _{1}^{\prime }(\kappa )=
\kappa /(5H_{0}+\kappa (2+5\nu_{0}))$.
\end{proposition}

\begin{remark} \label{rem:2} 
It is shown in Remark \ref{rem:6} that the Gaussian process
indexed by $\mathcal{G}\left( \kappa \right) $ 
\begin{equation*}
t^{\kappa }\mathbb{G}_{\left( 0,T\right) }\left( h_{t,x}\right) =t^{\kappa
}G\left( t,x\right) ,\left( t,x\right) \in \left[ 0,T\right] \times
\mathbb{R},
\end{equation*}
has a version that is uniformly continuous with bounded trajectories.
Therefore Proposition \ref{prop:2} implies that for any $\kappa >0$ the
weighted empirical process $t^{\kappa }\alpha _{n}\left( h_{t,x}\right)
=t^{\kappa }v_{n}\left( t,x\right) ,$ $\left( t,x\right) \in
\left[ 0,T\right] \times \mathbb{R}$, converges weakly to 
$t^{\kappa }G\left( t,x\right)$.
Recall, as pointed out in Remark \ref{rem:1}, weak convergence
fails if $\kappa$ is chosen to be zero.
\end{remark}

Propositions \ref{prop:1} and \ref{prop:2} lead to the following two strong
approximation theorems.

\begin{theorem}
\label{th:1} As long as $1\geq $ $\ \gamma =\gamma _{n}>0$ is constant,
under the assumptions and notation of Proposition \ref{prop:1} for all
$1/\left( 2\tau _{1}(0)\right) <\alpha <1/\tau _{1}(0)$ and $\xi >1$ there
exist a $\rho \left( \alpha ,\xi \right) >0$, a sequence of
i.i.d.~$B_{1}^{\left( H\right) },B_{2}^{(H)},\ldots ,$ and a sequence of 
independent
copies $\mathbb{G}_{\left( \gamma ,T\right) }^{\left( 1\right) },
\mathbb{G}_{\left( \gamma ,T\right) }^{\left( 2\right) },\ldots $, of 
$\mathbb{G}_{\left( \gamma ,T\right) }$ sitting on the same probability space 
such that 
\begin{equation*}
P\bigg\{ \max_{1\leq m\leq n}
\Big\Vert \sqrt{m}\alpha _{m} - \sum_{i=1}^{m} \mathbb{G}_{(\gamma, T)}^{(i)} 
\Big\Vert_{\mathcal{F}_{\left( \gamma ,T\right) }} > \rho ( \alpha ,\xi)
n^{1/2-\tau (\alpha) } \left( \log n\right)^{\tau_2}\bigg\}
\leq n^{-\xi }
\end{equation*}%
and 
\begin{equation*}
\max_{1\leq m\leq n} \Big\Vert \sqrt{m}\alpha_m- 
\sum_{i=1}^{m}\mathbb{G}_{(\gamma, T)}^{(i)} 
\Big\Vert _{\mathcal{F}_{(\gamma ,T)}}
=O\left( n^{1/2-\tau ( \alpha )} \left( \log n\right) ^{\tau _{2}}\right), \text{ a.s.,}
\end{equation*}%
where $\tau \left( \alpha \right) =
\left( \alpha \tau _{1}(0)-1/2\right)/(1+\alpha )>0.$
\end{theorem}

\begin{theorem}
\label{th:2} Under the assumptions and notation of Proposition \ref{prop:2}
for any $\kappa >0$, for all $1/\left( 2\tau _{1}^{\prime }\right) <\alpha
<1/\tau _{1}^{\prime }$, and $\xi >1$ there exist a $\rho ^{\prime }\left(
\alpha ,\xi \right) >0$, a sequence of i.i.d.~$B_{1}^{\left( H\right)
},B_{2}^{(H)},\ldots ,$ and a sequence of independent copies $\mathbb{G}%
_{\left( 0,T\right) }^{\left( 1\right) },\mathbb{G}_{\left( 0,T\right)
}^{\left( 2\right) },$ $\dots ,$ of $\mathbb{G}_{\left( 0,T\right) }$
sitting on the same probability space such that 
\begin{equation*}
P\left\{ \max_{1\leq m\leq n}
\Big\Vert \sqrt{m}\alpha_{m} - \sum_{i=1}^{m} \mathbb{G}_{(0,T)}^{(i)}
\Big\Vert_{\mathcal{G}(\kappa)} > 
\rho^{\prime}(\alpha ,\xi) n^{1/2-\tau^{\prime }(\alpha)} \left( \log 
n\right)^{\tau_2} \right\} \leq n^{-\xi }
\end{equation*}%
and 
\begin{equation}
\max_{1\leq m\leq n}
\Big\Vert \sqrt{m}\alpha_m - \sum_{i=1}^{m}\mathbb{G}_{(0,T)}^{(i)} 
\Big\Vert_{\mathcal{G}(\kappa) }
=O\left( n^{1/2-\tau ^{\prime }(\alpha) }\left(\log n\right)^{\tau_2}\right),
\text{ a.s.,}  \label{t3}
\end{equation}%
where $\tau ^{\prime }\left( \alpha \right) =
\left( \alpha \tau _{1}^{\prime}-1/2\right) /(1+\alpha )>0.$
\end{theorem}

\begin{remark}
\label{rem:3} Theorems \ref{th:1} and \ref{th:2} are strong approximations,
meaning that strong limit theorems can be inferred for the approximated
empirical process $\alpha_{n}$ from those that may hold for the sequence of
approximating Gaussian processes as long as the almost sure rate of strong
approximation is close enough. This is illustrated in Section \ref{sFLIL}.
\end{remark}

\subsection{Comments on the proofs of Theorems 1 and 2}
\label{sketch}

The proofs of Theorems \ref{th:1} and \ref{th:2} follow from
Propositions \ref{prop:1} and \ref{prop:2} (after some obvious notation
translations) exactly as Theorem 1 in \cite{BerthetMason} follows from their
Proposition 1, where a scheme described on pages 236--238 of Philipp
\cite{Phil} is closely followed. (Note that in \cite{BerthetMason}
\textquotedblleft $C\rho \left( \alpha ,\gamma \right) $\textquotedblright\
should be \textquotedblleft $\rho \left( \alpha ,\gamma \right) $%
\textquotedblright .) The essential ingredients are the maximal Inequalities
1A and 2A. Subsection \ref{four}.

\subsection{Applications to FLIL}
\label{sFLIL}

Theorem \ref{th:1} obviously implies that for any fixed choice of $0<\gamma
\leq 1<T$ there exist on the same probability space an i.i.d.~sequence $%
B_{1}^{(H)},B_{2}^{(H)},\ldots ,$ of sample continuous fractional Brownian
motions on $\left[ \gamma ,T\right] $ with Hurst index $0<H<1$ and a
sequence of independent copies $\mathbb{G}_{\left( \gamma ,T\right)
}^{\left( 1\right) },\mathbb{G}_{\left( \gamma ,T\right) }^{\left( 2\right)
},$ $\dots ,$ of $\mathbb{G}_{\left( \gamma ,T\right) }$ such that 
\begin{equation}
\begin{split}
\max_{1\leq m\leq n}
\Big\Vert \sqrt{m}\alpha_m - \sum_{i=1}^{m}\mathbb{G}_{( \gamma ,T)}^{(i)}
\Big\Vert_{\mathcal{F}_{(\gamma ,T)}} 
& = \max_{1\leq m\leq n} \sup_{(t,x) \in \mathcal{T}(\gamma)}
\big\vert \sqrt{m} v_{m}(t,x) - \sum_{i=1}^{m} G_{i}(t,x) \big\vert \\
& =o\left( \sqrt{n\log \log n }\right), \text{ a.s.,}  \label{AS}
\end{split}
\end{equation}%
where $\mathbb{G}_{\left( \gamma ,T\right) }^{\left( i\right) }\left(
h_{t,x}\right) =:G_{i}\left( t,x\right) $, for $i\geq 1$. Noting by the
comment right after Remark \ref{rem:c}, we can consider that each $%
G_{i}\left( t,x\right) $ is w.p.~$1$ [with probability $1$] in the separable
Banach space consisting of continuous functions in the sup-norm on the
compact metric space $\mathcal{T}^{c}\left( \gamma \right) $, equipped with
the semimetric $\rho $, we can apply the theorem in LePage \cite{LeP} (see
also Corollary 2.2 of Arcones \cite{Arcones}) to conclude the following
FLIL, namely, the sequence of Gaussian processes defined on $\mathcal{T}%
^{c}\left( \gamma \right)$ 
\begin{equation*}
\left\{ \frac{\sum_{i=1}^{n}G_{i}\left( t,x\right) }{\sqrt{2n\log \log n}}%
:\left( t,x\right) \in \mathcal{T}^{c}\left( \gamma \right) \right\}
\end{equation*}%
is w.p.~$1$ relatively compact in $\ell _{\infty }\left( \mathcal{T}%
^{c}\left( \gamma \right) \right) $, (the space of bounded functions $%
\Upsilon$ on $\mathcal{T}^{c}\left( \gamma \right) $ equipped with supremum
norm $\left\Vert \Upsilon \right\Vert _{\ell _{\infty }\left( \mathcal{T}%
^{c}\left( \gamma \right) \right) }=\sup_{\varphi \in \ell _{\infty }\left( 
\mathcal{T}^{c}\left( \gamma \right) \right) }\left\vert \Upsilon \left(
\varphi \right) \right\vert $), and its limit set is the unit ball of the
reproducing kernel Hilbert space determined by the covariance function $%
E\left( G\left( s,x\right) G\left( t,y\right) \right) .$ Note that by
continuity of $G\left( t,x\right) $ and its covariance function, the same
statement holds with $\mathcal{T}^{c}\left( \gamma \right) $ replaced by $%
\mathcal{T}\left( \gamma \right)$. Therefore by (\ref{AS}) the same is true
for 
\begin{equation}
\left\{ \frac{v_{n}\left( t,x\right) }{\sqrt{2\log \log n}}:\left(
t,x\right) \in \mathcal{T}\left( \gamma \right) \right\} .  \label{LL2}
\end{equation}%
This result can also be inferred from the FLIL for the empirical process as
stated in Theorem 9 on p.~609 of Ledoux and Talagrand
\cite{LedeouxTalagrand2} using the fact pointed out above that $v_{n}$ converges
weakly to a bounded uniformly continuous centered Gaussian process
$G\left(t,x\right) $ indexed by 
$\left( t,x\right) \in \mathcal{T}\left( \gamma\right)$. In particular we get 
that
\begin{equation*}
\limsup_{n\rightarrow \infty }\frac{\left\Vert \alpha _{n}\right\Vert _{%
\mathcal{F}_{\left( \gamma ,T\right) }}}{\sqrt{2\log \log n}}%
=\limsup_{n\rightarrow \infty }\sup_{\left( t,x\right) \in \mathcal{T}\left(
\gamma \right) }\left\vert \frac{v_{n}\left( t,x\right) }{\sqrt{2\log \log n}%
}\right\vert =\sigma \left( \gamma ,T\right),\ \text{ a.s.}
\end{equation*}%
where 
\begin{align*}
\sigma ^{2}\left( \gamma ,T\right) & =\sup \left\{ E\left( \mathbb{G}%
_{\left( \gamma ,T\right) }^{2}\left( h_{t,x}\right) \right) :h_{t,x}\in 
\mathcal{F}_{\left( \gamma ,T\right) }\right\} \\
& =\sup \left\{ \mathop{Var}(h_{t,x}(B^{(H)})):\left( t,x\right) \in 
\mathcal{T}\left( \gamma \right) \right\} =\frac{1}{4}.
\end{align*}

In the same way, on the probability space of Theorem \ref{th:2}, for all
$0<\kappa <\infty$, 
\begin{equation}
\begin{split}
\max_{1\leq m\leq n}
\Big\Vert \sqrt{m}\alpha_m - \sum_{i=1}^{m} \mathbb{G}_{(0,T)}^{(i)}
\Big\Vert_{\mathcal{G}(\kappa)} 
& =\max_{1\leq m\leq n} \sup_{(t,x) \in \mathcal{T} (0)}
\big\vert \sqrt{m} t^{\kappa }v_{m}(t,x) - \sum_{i=1}^{m} t^{\kappa }
G_{i}(t,x) \big\vert \\
& =o\left( \sqrt{n\log \log n}\right) ,\ \text{a.s.,}  \label{ASK}
\end{split}%
\end{equation}%
where $t^{\kappa }\mathbb{G}_{\left( 0,T\right) }^{\left( i\right) }
\left(h_{t,x}\right) =:t^{\kappa }G_{i}\left( t,x\right)$, for $i\geq 1$. We
point out in Remark \ref{rem:6} below that the process
$G_{\kappa }\left(t,x\right) :=t^{\kappa }G\left( t,x\right) $ has a version 
that is bounded
and uniformly continuous on $\mathcal{T}\left( 0\right) =[0,T]\times
\mathbb{R}$ with respect to the semimetric 
\begin{equation}
\rho _{\kappa }\left( (s,x\right) ,(t,y))=\sqrt{E\left( s^{\kappa }G\left(
s,x\right) -t^{\kappa }G\left( t,y\right) \right) ^{2}}.  \label{mk}
\end{equation}
\textit{From now on we assume that $G_{\kappa }\left( t,x\right)$
is such a version}. Denote by $\mathcal{T}^{c}\left( 0\right) $ the
completion of $\mathcal{T}\left( 0\right) $ in the topology induced by the
semimetric $\rho _{\kappa }$ from which we get by arguing as above and
applying the LePage theorem that%
\begin{equation*}
\left\{ \frac{\sum_{i=1}^{n}t^{\kappa }G_{i}\left( t,x\right) }{\sqrt{2n\log
\log n}}:\left( t,x\right) \in \mathcal{T}^{c}\left( 0\right) \right\}
\end{equation*}%
is, w.p.~$1$, relatively compact in
$\ell _{\infty }\left( \mathcal{T}^{c}\left( 0\right) \right)$ and its limit 
set is the unit ball of the
reproducing kernel Hilbert space determined by the covariance function
$E\left( s^{\kappa }t^{\kappa }G\left( s,x\right) G\left( t,y\right) \right)$,
$\left( s,x\right) \in $ $\mathcal{T}^{c}\left( 0\right)$. Note that by
uniform continuity of $G_{\kappa }\left( t,x\right) =
t^{\kappa }G\left(t,x\right) $ and its covariance function, the same statement 
holds with $\mathcal{T}^{c}\left( 0\right) $ replaced by
$\mathcal{T}\left( 0\right)$.
Therefore by (\ref{ASK}) the same is true for the sequence of processes 
\begin{equation}
\left\{ \frac{t^{\kappa }v_{n}\left( t,x\right) }{\sqrt{2\log \log n}}
:\left( t,x\right) \in \mathcal{T}\left( 0\right) \right\}.  \label{comL}
\end{equation}
This implies that 
\begin{equation}
\limsup_{n\rightarrow \infty }\sup_{\left( t,x\right) \in \mathcal{T}\left(
0\right) }\left\vert \frac{t^{\kappa }v_{n}\left( t,x\right) }{\sqrt{2\log
\log n}}\right\vert =\sigma _{\kappa }\left( T\right) ,\ \text{ a.s.}
\label{sig}
\end{equation}%
where 
\begin{equation*}
\begin{split}
\sigma _{\kappa }^{2}\left( T\right) & =\sup \left\{ E\left( \mathbb{G}%
_{\left( 0,T\right) }^{2}\left( t^{\kappa }h_{t,x}\right) \right) :t^{\kappa
}h_{t,x}\in \mathcal{G}\left( \kappa \right) \right\} \\
& =\sup \left\{ \mathop{Var}(t^{\kappa }h_{t,x}(B^{\left( H\right)
})):\left( t,x\right) \in \mathcal{T}\left( 0\right) \right\} =
\frac{T^{2\kappa }}{4}.
\end{split}
\end{equation*}
FLILs are by no means the only strong limit theorems for $\alpha _{n}$ that
can be derived from Theorems \ref{th:1} and \ref{th:2}. For instance, one
could consider Chung-type LILs.

\section{ Proofs of Propositions \protect\ref{prop:1} and
\protect\ref{prop:2}}

\label{s4}

Before we can prove Propositions \ref{prop:1} and \ref{prop:2} we must first
establish Proposition \ref{propPrime} below, which is a version of the
coupling given in Proposition \ref{prop:1} that holds on an appropriate
class of functions $\mathcal{F}_{n}$. To do so we must first define this
class of functions, derive an entropy bound for it and choose a good grid.
Our entropy bound will allow us to fill in the interstices of the empirical
and Gaussian processes constructed on $\mathcal{F}_{n}$ in Proposition
\ref{propPrime} by processes defined on all of $\mathcal{F}_{(\gamma _{n},T)}$
in such a way as to get useful rates of coupling. The proofs of the
bracketing bounds given in Subsection \ref{BBs} form the most technical part
of this paper.

\subsection{A useful class of functions}

To ease the notation from now on we suppress the Hurst index $H$. As above,
let $B\left( s\right) = B^{(H)}(s)$, $s\geq 0$, denote a sample continuous
fractional Brownian motion with Hurst index $0<H<1$. We have 
\begin{equation*}
E\left( B\left( t\right) -B\left( s\right) \right) ^{2}=\left\vert
t-s\right\vert ^{2H}.
\end{equation*}%
Note that for all 
$\left( s,x\right) ,\left( t,y\right) \in \mathcal{T}\left(\gamma \right)$, 
\begin{equation*}
\begin{split}
\rho ^{2}\left( (s,x\right) ,(t,y))& =E\left( 1\left\{ B(s)\leq x\right\}
-F(s,x)-\left( 1\left\{ B\left( t\right) \leq y\right\} -F\left( t,y\right)
\right) \right) ^{2} \\
& \leq E\left( 1\left\{ B\left( s\right) \leq x\right\} -1\left\{ B\left(
t\right) \leq y\right\} \right) ^{2}=d_{P}^{2}\left( h_{s,x},h_{t,y}\right) .
\end{split}
\end{equation*}

For the modulus of continuity of a sample continuous fractional Brownian
motion $B$, with Hurst index $H$, Wang (\cite{Wang}, Corollary 1.1) proved
that 
\begin{equation}
\lim_{h\downarrow 0}\sup_{t\in (0,1-h)}\frac{|B(t+h)-B(t)|}
{h^{H}\sqrt{2\log h^{-1}}}=1,\ \text{a.s.}  \label{W}
\end{equation}%
Recall the definition of $f_{H}$ in (\ref{FH}). For any $K\geq 1$ denote the
class of continuous real-valued functions on $\left[ 0,T\right] $, 
\begin{equation}
\mathcal{C}\left( K\right) =\left\{ g:\ \left\vert g(s)-g(t)\right\vert \leq
Kf_{H}(|s-t|),0\leq s,t\leq T\right\} .  \label{cck}
\end{equation}%
One readily checks that $\mathcal{C}\left( K\right)$ is closed in 
$\mathcal{C}\left[ 0,T\right]$. For any 
$\left( t,x\right) \in \mathcal{T}\left( \gamma \right)$ let 
$h_{t,x}^{\left(K\right) }$ denote the function of
$g\in \mathcal{C}\left[ 0,T\right] \rightarrow \left\{ 0,1\right\} $
defined by%
\begin{equation*}
h_{t,x}^{\left( K\right) }\left( g\right) =1\left\{ g\left( t\right) \leq
x,g\in \mathcal{C}\left( K\right) \right\} .
\end{equation*}%
The following class of functions will play an essential role in our proof: 
\begin{equation}
\mathcal{F}\left( K,\gamma \right) :=\left\{ h_{t,x}^{\left( K\right)
}:\left( t,x\right) \in \mathcal{T}\left( \gamma \right) \right\} .
\label{FK}
\end{equation}%
It is shown in the Appendix that these classes are \textit{pointwise
measurable}, which allows us to take supremums over these classes without
the need to worry about measurability problems.

\subsection{Bracketing}

We shall use the notion of bracketing. Let $\mathcal{G}$ be a class of
measurable real-valued functions defined on a measurable space $(S,\mathcal{S%
})$. A way to measure the size of a class $\mathcal{G}$ is to use $L_{2}(P)$%
-brackets. Let $l$ and $v$ be measurable real-valued functions on $\left( S,%
\mathcal{S}\right) $ such that $l\leq v$ and $d_{P}(l,v)<u,$ $u>0$, where 
\begin{equation*}
d_{P}(l,v)=\sqrt{E_{P}\left( l\left( \xi \right) -v\left( \xi \right)
\right) ^{2}}
\end{equation*}%
and $\xi $ is a random variable taking values in $S$ defined on a
probability space $\left( \Omega ,\mathcal{A},P\right) $. The pair of
functions $l$, $v$ form an $u$-bracket $\left[ l,v\right] $ consisting of
all the functions $f\in \mathcal{G}$ such that $l\leq f\leq v$. Let $%
N_{[\;]}(u,\mathcal{G},d_{P})$ be the minimum number of $u$-brackets needed
to cover $\mathcal{G}$.

\subsection{A useful bracketing bound}

\label{BBs}Our next aim is to bound the bracketing number
$N_{[\;]}(u,\mathcal{F}(K,\gamma),d_{P})$, where $P$ is the measure induced on 
the Borelsets of $\mathcal{C}\left[ 0,T\right] $, by $B$, with
$d_{P}^{2}\left(l,v\right) =
E\left( l\left( B\right) -v\left( B\right) \right)^{2}$.
\smallskip

We shall prove the following entropy bound:\smallskip

\noindent\textbf{Entropy Bound I} \ For some constant $C_{T}$
(depending on $T$ and $H$), for $u\in\left( 0,1/e\right)$,
$\gamma\in\left( 0,1/e\right) $
and $K\geq e$,
\begin{equation}
N_{[\;]}(u,\mathcal{F}(K,\gamma),d_{P})\leq C_{T}K^{1/H}u^{-2(1+1/H)}
\sqrt{\log u^{-1}}\gamma^{-(1+H)}\left( \log\left( \frac{K}{u\gamma}\right)
\right) ^{\frac{1}{H}}.  \label{bn1}
\end{equation}

\noindent\textit{Proof} 
Choose $\gamma=t_{0}<t_{1}<\ldots<t_{k}=T$, such that 
\begin{equation}
Kf_{H}\left( t_{i}-t_{i-1}\right) \leq 1, \text{ for }0\leq i\leq k,
\label{t}
\end{equation}
and $x_{-m}<x_{-m+1}<\ldots<x_{-1}<x_{0}=0<x_{1}<\ldots<x_{m}$, with
$0=y_{0}<y_{1}<\ldots<y_{m}$, $x_{\pm i}=\pm y_{i}$, $i=0,1,\ldots,m$, such
that 
\begin{equation}
x_{m}\geq2T^{H}.  \label{xm}
\end{equation}
Also put $x_{-(m+1)}=-\infty$, $x_{m+1}=\infty$.

Consider the upper and lower functions: for $g\in\mathcal{C}[0,T]$ 
\begin{equation*}
v_{i,j}(g)=1 \left\{ g(t_{i-1}) \leq x_{j} + Kf_{H}\left(
t_{i}-t_{i-1}\right) , g\in\mathcal{C}(K) \right\}
\end{equation*}
and 
\begin{equation*}
l_{i,j}(g)=1 \left\{ g(t_{i-1})\leq x_{j-1} - Kf_{H}\left(
t_{i}-t_{i-1}\right) , g\in\mathcal{C}(K) \right\} ,
\end{equation*}
for $i=1,2,\ldots,k$, $j=-m,\ldots,m,m+1$. Note that $v_{i,m+1}(g)=1\{g\in 
\mathcal{C}(K)\}$, and $l_{i,-m}(g)=0$.

First we show that these functions define a covering. Select any
$t_{i-1}<t\leq t_{i}$ (in the case $i=1$ we allow $t_{0}=t$) and
$x_{j-1}<x\leq x_{j}$, for $i=1,\dots,k,$ $j=-m+1,\dots,m.$ Since for any
$g\in\mathcal{C}\left( K\right) $ 
\begin{equation*}
g\left( t_{i-1}\right) -Kf_{H}(t_{i}-t_{i-1})\leq g\left( t\right) \leq
g\left( t_{i-1}\right) +Kf_{H}\left( t_{i}-t_{i-1}\right)
\end{equation*}
we see that for all $g\in\mathcal{C}\left( K\right) $, $l_{i,j}\left(
g\right) \leq h_{t,x}^{\left( K\right) }\left( g\right) \leq v_{i,j}\left(
g\right)$. \smallskip

Next, for $-\infty<x\leq x_{-m}$ and any $t_{i-1}<t\leq t_{i}$, $%
0=l_{i,-m}\left( g\right) \leq h_{t,x}^{\left( K\right) }\left( g\right)
\leq v_{i,-m}\left( g\right) $, and for $x_{m}<x<\infty$ and any $t_{i-1} <
t \leq t_{i},$ $l_{i,m+1}\left( g\right) \leq h_{t,x}^{\left( K\right)
}\left( g\right) \leq v_{i,m+1}\left( g\right) =1\{ g \in\mathcal{C}(K) \}$.

Clearly for $-m+1\leq j\leq m$ we get 
\begin{equation}
\begin{split}
& d_{P}^{2} (l_{i,j}, v_{i,j}) = E\left( v_{i,j}(B)-l_{i,j}(B)\right) ^{2} \\
& = P \Big\{ B(t_{i-1}) \in(x_{j-1} - Kf_{H}(t_{i}-t_{i-1}), x_{j} + K
f_{H}(t_{i}-t_{i-1})], B \in\mathcal{C}(K) \Big\} \\
& \leq P\left\{ B(t_{i-1})
\in(x_{j-1}-Kf_{H}(t_{i}-t_{i-1}),x_{j}+Kf_{H}(t_{i}-t_{i-1})]\right\} \\
& =\Phi\left( \frac{x_{j}+Kf_{H}(t_{i}-t_{i-1})}{t_{i-1}^{H}}\right)
-\Phi\left( \frac{x_{j-1}-Kf_{H}(t_{i}-t_{i-1})}{t_{i-1}^{H}}\right) .
\end{split}
\label{ph}
\end{equation}
So that for $-m+1\leq j\leq m$ we have 
\begin{equation*}
d_{P}^{2}(l_{i,j},v_{i,j})\leq\frac{1}{\sqrt{2\pi}}\left(
x_{j}-x_{j-1}+2Kf_{H}(t_{i}-t_{i-1})\right) t_{i-1}^{-H}.
\end{equation*}
Inequality (\ref{ph}) is also valid for $j=-m$ and $j=m+1$, namely 
\begin{equation*}
d_{P}^{2}(l_{i,-m},v_{i,-m})=d_{P}^{2}(l_{i,m+1},v_{i,m+1})\leq1-\Phi\left( 
\frac{x_{m}-Kf_{H}(t_{i}-t_{i-1})}{t_{i-1}^{H}}\right) .
\end{equation*}
Now by $t_{i-1}^{H}\leq T^{H}$, $2T^{H}\geq2$, (\ref{t}) and (\ref{xm}) we
have%
\begin{equation*}
\frac{x_{m}-Kf_{H}(t_{i}-t_{i-1})}{t_{i-1}^{H}}=\frac{%
2x_{m}-2Kf_{H}(t_{i}-t_{i-1})}{2t_{i-1}^{H}}\geq\frac{x_{m}}{2T^{H}},
\end{equation*}
which when combined with the standard normal tail bound holding for $z>0$, $%
P\left\{ Z\geq z\right\} \leq\frac{1}{z\sqrt{2\pi}}\exp(-z^{2}/2),$ gives 
\begin{equation*}
1-\Phi\left( \frac{x_{m}-Kf_{H}(t_{i}-t_{i-1})}{t_{i-1}^{H}}\right)
\leq1-\Phi\left( \frac{x_{m}}{2T^{H}}\right) \leq\frac{1}{\sqrt{2\pi}}\frac{%
2T^{H}}{x_{m}}e^{-\frac{x_{m}^{2}}{8T^{2H}}}.
\end{equation*}
From this we see that for $u \in(0,e^{-1})$, the choice $x_{m}\geq4T^{H}%
\sqrt{\log u^{-1}}$ ensures that 
\begin{equation*}
d_{P}^{2}(l_{i,-m},v_{i,-m})=d_{P}^{2}(l_{i,m+1},v_{i,m+1})\leq u^{2}.
\end{equation*}
Thus to construct our $u$-covering, it suffices to appropriately partition
the intervals 
\begin{equation*}
\left[ -4T^{H}\sqrt{\log u^{-1}},4T^{H}\sqrt{\log u^{-1}}\right] \quad 
\text{and } \ [\gamma,T],
\end{equation*}
so that $x_{m}\geq4T^{H}\sqrt{\log u^{-1}}$, $t_{i}-t_{i-1}$ satisfies
(\ref{t}), and for $0\leq i\leq k$ and $-m+1\leq j\leq m$,
$d_{P}^{2}(l_{i,j},v_{i,j})\leq u^{2}.$ \smallskip

Set 
\begin{equation}
\Delta(\gamma,u)=\sqrt{\frac{\pi}{2}}\gamma^{H}u^{2}\text{ and }\Gamma
(\gamma,u)=\left( \sqrt{\frac{\pi}{8}}\right) ^{1/H}\frac{K^{-1/H}\gamma
u^{2/H}}{\left[ \log\left( K^{1/H}\gamma^{-1}u^{-2/H}\right) \right] ^{1/H}}.
\label{del2}
\end{equation}

Let $\left\lceil x\right\rceil $ denote here and elsewhere the smallest
integer greater than or equal to $x$. Putting 
\begin{equation*}
m(\gamma,u)=\left\lceil \frac{4T^{H}\sqrt{\log u^{-1}}}{\Delta(\gamma ,u)}%
\right\rceil =:m\ \text{ and }k(\gamma,u)=\left\lceil \frac{T-\gamma }{%
\Gamma(\gamma,u)}\right\rceil =:k,
\end{equation*}
straightforward computations show that for the choice 
\begin{equation*}
y_{i}=i\Delta(\gamma,u),i=0,1,\ldots,m, \ t_{j}=\gamma+j\Gamma
(\gamma,u),j=0,1,\ldots,k-1\text{ and }t_{k}=T,
\end{equation*}
by (\ref{del2}) we have for $-m+1\leq j\leq m$ 
\begin{equation*}
d_{P}^{2}(l_{i,j},v_{i,j}) \leq u^{2}.
\end{equation*}
We also see that $y_{m}=x_{m}=4T^{H}\sqrt{\log u^{-1}}\geq2T^{H}$, and by
(\ref{del2}) for $0\leq i\leq k$ 
\begin{equation*}
Kf_{H}\left( t_{i}-t_{i-1}\right) \leq Kf_{H}(\Gamma(\gamma,u))\leq \frac{%
\gamma^{H}u^{2}\sqrt{2\pi}}{4}\leq\gamma^{H}u^{2}\leq1.
\end{equation*}
Thus (\ref{t}) and (\ref{xm}) hold. Hence this choice of $t_{i}$ and $x_{j}$
corresponds to a $u$-covering of $\mathcal{F}(K,\gamma)$. So we have proved
the following entropy bound: for $u\in(0,e^{-1}),\gamma\in(0,e^{-1})$ and $%
K\geq e$ 
\begin{equation*}
N_{[\;]}(u,\mathcal{F}(K,\gamma),d_{P}) \leq\left( k(\gamma,u)+1\right)
(2m(\gamma,u)+2),
\end{equation*}
thus (\ref{bn1}) holds for some constant $C_{T}$ (depending on $T$ and $H$).
\hfill$\square$ \smallskip

It will often be convenient to use the following weaker entropy bound, which
follows easily from (\ref{bn1}). \medskip

\noindent\textbf{Entropy Bound II} For some constant $C_{T}^{\prime}$
(depending on $T$ and $H$), for $u\in\left( 0,1\right] ,$ $\gamma\in\left(
0,1\right] $ and $K\geq e,$ 
\begin{equation}
N_{[\;]}(u,\mathcal{F}(K,\gamma),d_{P})\leq
C_{T}^{\prime}K^{2/H}u^{-3(1+1/H)}\gamma^{-(1+2/H)}.  \label{eq:weak-bn1}
\end{equation}
\smallskip

Set 
\begin{equation}
\mathcal{F}\left( K,\gamma,\varepsilon\right) =\left\{ \left( f,f^{\prime
}\right) \in\mathcal{F}\left( K,\gamma\right) \times\mathcal{F}\left(
K,\gamma\right) :d_{P}\left( f,f^{\prime}\right) <\varepsilon\right\}
\label{F}
\end{equation}
and%
\begin{equation}
\mathcal{G}\left( K,\gamma,\varepsilon\right) =\left\{ f-f^{\prime}:\left(
f,f^{\prime}\right) \in\mathcal{F}\left( K,\gamma,\varepsilon\right)
\right\} ,  \label{G}
\end{equation}
that is, $\mathcal{F}\left( K,\gamma,\varepsilon\right) $ and $\mathcal{G}%
\left( K,\gamma,\varepsilon\right) $ are the classes of functions on $%
\mathcal{C}\left[ 0,T\right] $, indexed by $\gamma\leq s,t\leq T$, $%
-\infty<x,y<\infty,$ defined for $g\in\mathcal{C}\left[ 0,T\right] $ by 
\begin{equation*}
\left( h_{s,x}^{\left( K\right) }\left( g\right) ,h_{t,y}^{\left( K\right)
}\left( g\right) \right) =\left( 1\left\{ g\left( s\right) \leq x,g\in%
\mathcal{C}\left( K\right) \right\} ,1\left\{ g\left( t\right) \leq y,g\in%
\mathcal{C}\left( K\right) \right\} \right)
\end{equation*}
and%
\begin{equation*}
h_{s,x}^{\left( K\right) }\left( g\right) -h_{t,y}^{\left( K\right) }\left(
g\right) ,
\end{equation*}
respectively, and satisfying 
\begin{equation*}
d_{P}\left( h_{s,x}^{\left( K\right) },h_{t,y}^{\left( K\right) }\right) =%
\sqrt{E\left( h_{s,x}^{\left( K\right) }\left( B\right) -h_{t,y}^{\left(
K\right) }\left( B\right) \right) ^{2}}<\varepsilon.
\end{equation*}
\medskip

\noindent We find that independently of $\varepsilon$ 
\begin{equation}
N_{[\;]}(u,\mathcal{G}\left( K,\gamma,\varepsilon\right) ,d_{P})\leq\left(
N_{[\;]}(u/2,\mathcal{F}\left( K,\gamma\right) ,d_{P})\right) ^{2}.
\label{NN}
\end{equation}

\subsection{Proof of Proposition \protect\ref{prop:1}}

\label{prp1}For any $c>0$, $n>e$ and $0<\gamma \leq 1<T$ denote the class of
real-valued functions on $\left[ 0,T\right] $, 
\begin{equation}
\mathcal{C}_{n}:=\mathcal{C}(\sqrt{c\log n})=\left\{ g:\,\left\vert g\left(
s\right) -g\left( t\right) \right\vert \leq \sqrt{c\log n}f_{H}(\left\vert
s-t\right\vert ),\,0\leq s,t\leq T\right\} ,  \label{cn}
\end{equation}%
and let $\mathcal{C}_{\infty }$ be as in (\ref{CI}). Notice that by (\ref{W}),
$P\left\{ B \in \mathcal{C}_{\infty }\right\} =1$. Define the class of
functions $\mathcal{C}\left[ 0,T\right] \rightarrow \mathbb{R}$ indexed by
$\left[ \gamma _{n},T\right] \times \mathbb{R} = \mathcal{T}\left(
\gamma_{n}\right)$ 
\begin{equation*}
\mathcal{F}_{n}=\left\{ h_{t,x}^{\left( \sqrt{c\log n}\right) }\left(
g\right) =1\left\{ g\left( t\right) \leq x,g\in \mathcal{C}_{n}\right\}
:\left( t,x\right) \in \mathcal{T}\left( \gamma _{n}\right) \right\} .
\end{equation*}%
To simplify our previous notation we shall write here 
\begin{equation}
h_{t,x}^{(n)}\left( g\right) = h_{t,x}^{\left( \sqrt{c\log n}\right)
}\left(g\right).  \label{notation1}
\end{equation}%
For $h_{t,x}^{(n)}\in \mathcal{F}_{n}$ let 
\begin{equation*}
\alpha _{n}\left(h_{t,x}^{(n)}\right) =n^{-1/2}\sum_{i=1}^{n}\left( 1\left\{
B_{i}(t)\leq x,B_{i}\in \mathcal{C}_{n}\right\} -P\left\{ B(t)\leq x,B\in 
\mathcal{C}_{n}\right\} \right) .
\end{equation*}%
Notice that for each $\left( t,x\right) \in \mathcal{T}\left( \gamma
_{n}\right) $, when $B_{i}\in \mathcal{C}_{n}$, for $i=1,\dots ,n$, 
\begin{equation}
\begin{split}
\alpha _{n} \left( h_{t,x}^{(n)}\right) &= v_{n}\left( t,x\right) +\sqrt{n}
P\left\{ B\left( t\right) \leq x,B\notin \mathcal{C}_{n}\right\} \\
&= \alpha_{n}\left( h_{t,x}\right) + \sqrt{n}P\left\{ B \left( t\right) \leq
x,B\notin \mathcal{C}_{n}\right\}.  \label{note}
\end{split}%
\end{equation}%
Let $\mathbb{F}_{\left( \gamma _{n},T\right) }^{(n)}$ denote the mean zero
Gaussian process indexed by $\mathcal{F}_{n}$, having covariance function
defined for $h_{s,x}^{(n)},h_{t,y}^{(n)}\in \mathcal{F}_{n}$ by 
\begin{equation*}
\begin{split}
& E\left( \mathbb{F}_{\left( \gamma _{n},T\right) }^{(n)}\left(
h_{s,x}^{(n)}\right) \mathbb{F}_{\left( \gamma _{n},T\right) }^{(n)}\left(
h_{t,y}^{(n)}\right) \right) =P\left\{ B\left( s\right) \leq x,B\left(
t\right) \leq y,B\in \mathcal{C}_{n}\right\} \\
& \phantom{=}-P\left\{ B\left( s\right) \leq x,B\in \mathcal{C}_{n}\right\}
P\left\{ B\left( t\right) \leq y,B\in \mathcal{C}_{n}\right\} .
\end{split}%
\end{equation*}

We shall first establish the following auxiliary result.

\begin{proposition} \label{propPrime} 
As long as $1\geq $ $\ \gamma =\gamma_{n}>0$ satisfies (\ref{eta}),
for every $\vartheta >1$ there exists a 
$\eta \left( \vartheta \right) >0$ such that for each integer $n$ large enough 
one can construct on
the same probability space random vectors $B_{1}, \ldots,B_{n}$ i.i.d.~$B$
and a probabilistically equivalent version $\widetilde{\mathbb{F}}_{\left(
\gamma _{n},T\right) }^{(n)}$ of $\mathbb{F}_{\left( \gamma _{n},T\right)
}^{(n)}$ such that 
\begin{equation}
P\left\{ \left\Vert \alpha _{n}-\widetilde{\mathbb{F}}_{\left( \gamma
_{n},T\right) }^{(n)}\right\Vert _{\mathcal{F}_{n}}>\eta \left( \vartheta
\right) \left( \log n\right) ^{\tau _{2}}\left( n^{-1/2}\gamma
_{n}^{-5H_{0}/2}\right) ^{2/(2+5\nu _{0})}\right\} \leq n^{-\vartheta },
\label{p}
\end{equation}%
where $\tau _{2}$ is given in Proposition \ref{prop:1} and $H_{0}$ and
$\nu_{0}$ are defined as in (\ref{nuk}). Moreover, in particular, when
$\gamma_{n}=n^{-\eta }$, with $0<\eta <\frac{1}{5H_{0}}$ and is $\tau _{1}$ as 
in
Proposition \ref{prop:1}, 
\begin{equation}
P\left\{ \left\Vert \alpha _{n}-\widetilde{\mathbb{F}}_{\left( \gamma
_{n},T\right) }^{(n)}\right\Vert _{\mathcal{F}_{n}}>\eta \left( \vartheta
\right) n^{-\tau _{1}}\left( \log n\right) ^{\tau _{2}}\right\} \leq
n^{-\vartheta }.  \label{pp}
\end{equation}
\end{proposition}

\noindent \textit{Proof} 
Let $B$ be a sample continuous fractional Brownian
motion with Hurst index $0<H<1$ restricted to $\left[ 0,T\right] $ taking
values in the measurable space $(\mathcal{C}\left[ 0,T\right] ,\mathcal{B}%
\left[ 0,T\right] )$. As above $P$ denotes the probability measure induced
on the Borel sets of $\mathcal{C}\left[ 0,T\right] $ by $B$. Let $\mathcal{M}$
denote the real-valued measurable functions on the space
$(\mathcal{C}\left[ 0,T\right] ,\mathcal{B}\left[ 0,T\right])$. For any 
$\varepsilon >0$
we can choose a grid 
\begin{equation}
\mathcal{H}\left( \varepsilon \right) =\left\{ h_{k}:1\leq k\leq N\left(
\varepsilon \right) \right\}  \label{Hf}
\end{equation}%
of measurable functions $\mathcal{M}$ on 
$(\mathcal{C}\left[ 0,T\right],\mathcal{B}\left[ 0,T\right] )$ such that each 
$f\in \mathcal{F}_{n}$ is in
a ball $\{f\in \mathcal{M}:d_{P}(h_{k},f)<\varepsilon \}$ around some
$h_{k}$, $1\leq k\leq N\left( \varepsilon \right)$, where 
\begin{equation*}
d_{P}(h_{k},f)=
\sqrt{E\left( h_{k}\left( B\right) -f\left( B\right) \right)^{2}}.
\end{equation*}%
The choice 
\begin{equation}
N(\varepsilon )=N_{\left[ \;\right] }(\varepsilon /2,\mathcal{F}_{n},d_{P})
\label{Ne}
\end{equation}%
permits us to select such $h_{k}$ $\in \mathcal{F}_{n}$. Recalling the
previous notation (\ref{F}) and (\ref{G}), set 
\begin{equation}
\mathcal{F}_{n}\left( \varepsilon \right) =\mathcal{F}\left( \sqrt{c\log n}%
,\gamma _{n},\varepsilon \right) \text{ and }\mathcal{G}_{n}\left(
\varepsilon \right) =\mathcal{G}\left( \sqrt{c\log n},\gamma
_{n},\varepsilon \right) .  \label{notation}
\end{equation}%
Fix $n\geq 1$. Let $B_{1},\ldots ,B_{n}$ be i.i.d.~$B$, and
$\epsilon_{1},\ldots ,\epsilon_{n}$ be independent Rademacher random
variables mutually independent of $B_{1},\ldots,B_{n}$. Write for 
$\varepsilon >0$, 
\begin{equation}
\begin{split}
\mu _{n}^{S}\left( \varepsilon \right) & =E\left\{ \sup_{\left( f,f^{\prime
}\right) \in \mathcal{F}_{n}\left( \varepsilon \right) }\left\vert
n^{-1/2}\sum_{i=1}^{n}\epsilon _{i}\left( f-f^{\prime }\right) \left(
B_{i}\right) \right\vert \right\} \\
& =E\left\{ \sup_{f-f^{\prime }\in \mathcal{G}_{n}\left( \varepsilon \right)
}\left\vert n^{-1/2}\sum_{i=1}^{n}\epsilon _{i}\left( f-f^{\prime }\right)
\left( B_{i}\right) \right\vert \right\},
\end{split}
\label{H10}
\end{equation}%
and 
\begin{equation}
\mu _{n}^{G}\left( \varepsilon \right) =E\left\{ \sup_{\left( f,f^{\prime
}\right) \in \mathcal{F}_{n}\left( \varepsilon \right) }\left\vert \mathbb{F}%
_{\left( \gamma _{n},T\right) }^{(n)}(f)-\mathbb{F}_{\left( \gamma
_{n},T\right) }^{(n)}(f^{\prime })\right\vert \right\} .  \label{H101}
\end{equation}

\begin{lemma}
\label{lem:1} Given $\varepsilon >0$, $\delta >0$, $t>0$ and $n\geq 1$ large
enough, there exist a probability space $\left(\Omega ,\mathcal{A},P\right)$
on which sit $B_{1},\ldots,B_{n}$ i.i.d.~$B$ and a probabilistically
equivalent version $\widetilde{\mathbb{F}}_{\left( \gamma_{n},T\right)
}^{\left( n\right)}$ of the Gaussian process $\mathbb{F}_{\left( \gamma
_{n},T\right) }^{(n)}$ indexed by $\mathcal{F}_{n}$ such that for suitable
positive constants $C_{1},C_{2},A$, $A_{1}$ and $A_{5}$ with $A_{5} \leq 1/2$%
, independent of $\varepsilon >0$, $\delta >0$, $t>0$ and $n\geq 1$, we have 
\begin{equation}
\begin{split}
& P\left\{ \left\Vert \alpha _{n}-\widetilde{\mathbb{F}}_{\left( \gamma
_{n},T\right) }^{(n)}\right\Vert _{\mathcal{F}_{n}}>A\mu _{n}^{S}\left(
\varepsilon \right) +\mu _{n}^{G}\left( \varepsilon \right) +\delta +\left(
A+1\right) t\right\} \\
& \leq C_{1}N\left( \varepsilon \right) ^{2}\exp \left( -\frac{C_{2}\sqrt{n}%
\ \delta }{\left( N\left( \varepsilon \right) \right) ^{5/2}}\right) +2\exp
\left( -A_{1}\sqrt{n}\ t\right) +4\exp \left( -\frac{A_{5}t^{2}}{\varepsilon
^{2}}\right) .
\end{split}
\label{H9}
\end{equation}
\end{lemma}

\noindent \textit{Proof of Lemma \ref{lem:1}} Our proof applies the
procedure detailed in Section 5.1 in \cite{BerthetMason}. Given $\varepsilon
>0$ and $n\geq 1$, our aim is to construct a probability space $\left(
\Omega,\mathcal{A},P\right)$ on which sit $B_{1},\ldots ,B_{n}$ i.i.d.~$B$
and a version $\widetilde{\mathbb{F}}_{\left( \gamma _{n},T\right) }^{\left(
n\right) }$ of the Gaussian process $\mathbb{F}_{\left( \gamma _{n},T\right)
}^{(n)}$ indexed by $\mathcal{F}_{n}$ such that for $\mathcal{H}\left(
\varepsilon \right) $ and $\mathcal{F}_{n}\left( \varepsilon \right) $
defined as above and for all $A>0$, $\delta >0$ and $t>0$, 
\begin{equation}
\begin{split}
& P\left\{ \left\Vert \alpha _{n}-\widetilde{\mathbb{F}}_{\left( \gamma
_{n},T\right) }^{(n)}\right\Vert _{\mathcal{F}_{n}}>A\mu _{n}^{S}\left(
\varepsilon \right) +\mu _{n}^{G}\left( \varepsilon \right) +\delta +\left(
A+1\right) t\right\} \\
& \leq P\left\{ \max_{h\in \mathcal{H}\left( \varepsilon \right) }\left\vert
\alpha _{n}\left( h\right) -\widetilde{\mathbb{F}}_{\left( \gamma
_{n},T\right) }^{(n)}(h)\right\vert >\delta \right\} \\
& \phantom{\leq} \ +P\left\{ \sup_{\left( f,f^{\prime }\right) \in \mathcal{F%
}_{n}\left( \varepsilon \right) }\left\vert \alpha _{n} \left( f\right)
-\alpha_{n}\left( f^{\prime }\right) \right\vert >A \mu_{n}^{S}\left(
\varepsilon \right) +At\right\} \\
& \phantom{\leq} \ +P\left\{ \sup_{\left( f,f^{\prime }\right) \in \mathcal{F%
}_{n}\left( \varepsilon \right) }\left\vert \widetilde{\mathbb{F}}_{\left(
\gamma _{n},T\right) }^{(n)}(f)-\widetilde{\mathbb{F}}_{\left( \gamma
_{n},T\right) }^{(n)}(f^{\prime })\right\vert >t+\mu _{n}^{G}\left(
\varepsilon \right) \right\} \\
& =:P_{n}\left( \delta \right) +Q_{n}\left( t,\varepsilon \right) +%
\widetilde{Q}_{n}\left( t,\varepsilon \right),  \label{ine}
\end{split}%
\end{equation}
with all these probabilities simultaneously small for suitably chosen $A>0$, 
$\delta >0$ and $t>0$. Consider the $n$ i.i.d.~mean zero random vectors in $%
\mathbb{R}^{N\left( \varepsilon \right) }$, 
\begin{equation*}
Y_{i}:=\frac{1}{\sqrt{n}}\left( h_{1}\left( B_{i}\right) - Eh_{1}\left( B\right) ,
\dots , h_{N\left( \varepsilon \right) }\left(
B_{i}\right) - E h_{N\left( \varepsilon \right) }\left( B\right)
\right), \ 1\leq i\leq n.
\end{equation*}%
First note that by the definition of $h_{k}$ $\in \mathcal{F}_{n}$, we have 
\begin{equation*}
|Y_{i}|_{N\left( \varepsilon \right) }\leq \sqrt{\frac{N\left( \varepsilon
\right) }{n}},\mbox{ }1\leq i\leq n,
\end{equation*}%
where $\left\vert \cdot \right\vert _{N}$, $N\geq 1$, denotes the usual
Euclidean norm on $\mathbb{R}^{N}$. Therefore by the coupling inequality (\ref{coup}) we 
can enlarge the probability space on which (\ref{ine}) holds to
include $Z_{1},\dots ,Z_{n}$ i.i.d. 
\begin{equation*}
Z:=\left( Z^{1},\dots ,Z^{N\left( \varepsilon \right) }\right)
\end{equation*}%
mean zero Gaussian vectors such that 
\begin{equation}
P_{n}\left( \delta \right) \leq P\left\{ \left\vert
\sum_{i=1}^{n}\left( Y_{i}-Z_{i}\right) \right\vert _{N\left( \varepsilon
\right) }>\delta \right\} \leq C_{1}N\left( \varepsilon \right) ^{2}\exp
\left( -\frac{C_{2}\sqrt{n}\ \delta }{\left( N\left( \varepsilon \right)
\right) ^{5/2}}\right) ,  \label{bound1}
\end{equation}%
where $\mathop{Cov}(Z^{l},Z^{k})=\mathop{Cov}(Y^{l},Y^{k})=:\left\langle
h_{l},h_{k}\right\rangle $. Moreover by Lemma A1 of Berkes and Philipp this
space can be extended to include a\ probabilistically equivalent version $%
\widetilde{\mathbb{F}}_{\left( \gamma _{n},T\right) }^{(n)}$ of the Gaussian
process $\mathbb{F}_{\left( \gamma _{n},T\right) }^{(n)}$ indexed by $%
\mathcal{F}_{n}$ such that for $1\leq k\leq N\left( \varepsilon \right)$, 
\begin{equation*}
\widetilde{\mathbb{F}}_{\left( \gamma _{n},T\right)
}^{(n)}(h_{k})=\sum_{i=1}^{n}Z_{i}^{k}.
\end{equation*}%
The $P_{n}\left( \delta \right) $ in (\ref{ine}) is defined through this $%
\widetilde{\mathbb{F}}_{\left( \gamma _{n},T\right) }^{(n)}$. Notice that
the probability space on which $Y_{1},\dots ,Y_{n}$, $Z_{1},\dots ,Z_{n}$
and $\widetilde{\mathbb{F}}_{\left( \gamma _{n},T\right) }^{(n)}$ sit
depends on $n\geq 1$ and the choice of $\varepsilon >0$ and $\delta >0$.\smallskip

Observe that the class 
\begin{equation*}
\mathcal{G}_{n}\left( \varepsilon \right) =\left\{ f-f^{\prime }:\left(
f,f^{\prime }\right) \in \mathcal{F}_{n}\left( \varepsilon \right) \right\}
\end{equation*}%
satisfies 
\begin{equation*}
\sigma _{\mathcal{G}_{n}\left( \varepsilon \right) }^{2}=\sup_{\left(
f,f^{\prime }\right) \in \mathcal{F}_{n}\left( \varepsilon \right) }%
\mathop{Var}(f(B)-f^{\prime }(B))\leq \sup_{\left( f,f^{\prime }\right) \in 
\mathcal{F}_{n}\left( \varepsilon \right) }d_{P}^{2}\left( f,f^{\prime
}\right) \leq \varepsilon ^{2}.
\end{equation*}%
Thus with $A>0$ as in (\ref{tal}) we get by applying Talagrand's inequality, with $M=1$, 
\begin{equation}
Q_{n}\left( t,\varepsilon \right) =P\left\{ ||\alpha _{n}||_{%
\mathcal{G}_{n}\left( \varepsilon \right) }>A\left( \mu _{n}^{S}\left(
\varepsilon \right) +t\right) \right\} \leq 2\exp \left( -\frac{A_{1}t^{2}}{%
\varepsilon ^{2}}\right) +2\exp \left( -A_{1}\sqrt{n}t\right).
\label{bound2}
\end{equation}
Next, consider the separable centered Gaussian process $\mathbb{Z}_{\left(
f,f^{\prime }\right) }=\widetilde{\mathbb{F}}_{\left( \gamma _{n},T\right)
}^{(n)}(f)-\widetilde{\mathbb{F}}_{\left( \gamma _{n},T\right)
}^{(n)}(f^{\prime })$ indexed by $\mathcal{F}_{n}\left( \varepsilon \right)$. We have 
\begin{equation*}
\begin{split}
\sigma _{T}^{2}\left( \mathbb{Z}\right) & =\sup_{\left( f,f^{\prime }\right)
\in \mathcal{F}_{n} \left( \varepsilon \right) }  E \left( \left( 
\widetilde{\mathbb{F}}_{\left( \gamma _{n},T\right) }^{(n)}(f) -\widetilde{%
\mathbb{F}}_{\left( \gamma _{n},T\right) }^{(n)}(f^{\prime})
\right)^{2}\right) \\
& =\sup_{\left( f,f^{\prime }\right) \in \mathcal{F}_{n}\left( \varepsilon
\right) }\mathop{Var}(f(B)-f^{\prime }(B)) \leq \sup_{\left( f,f^{\prime
}\right) \in \mathcal{F}_{n}\left( \varepsilon \right) }d_{P}^{2}\left(
f,f^{\prime }\right) \leq \varepsilon^{2}.
\end{split}%
\end{equation*}%
Borell's inequality (\ref{bor}) now gives 
\begin{equation}
\widetilde{Q}_{n}\left( t,\varepsilon \right) = P \left\{
\sup_{\left( f,f^{\prime }\right) \in \mathcal{F}_{n}\left( \varepsilon
\right) }\left\vert \widetilde{\mathbb{F}}_{\left( \gamma _{n},T\right)
}^{(n)}(f)-\widetilde{\mathbb{F}}_{\left( \gamma _{n},T\right)
}^{(n)}(f^{\prime })\right\vert >t+\mu _{n}^{G}\left( \varepsilon \right)
\right\} \leq 2\exp \left( -\frac{t^{2}}{2\varepsilon ^{2}}\right) .
\label{bound3}
\end{equation}%
Putting (\ref{bound1}), (\ref{bound2}) and (\ref{bound3}) together we
obtain, for some positive constants $A$, $A_{1}$ and $A_{5}$ with $A_{5}\leq
1/2,$%
\begin{equation*}
\begin{split}
& P \left\{ \left\Vert \alpha _{n}-\widetilde{\mathbb{F}}_{\left(
\gamma _{n},T\right) }^{(n)}\right\Vert _{\mathcal{F}_{n}}>A\mu
_{n}^{S}\left( \varepsilon \right) +\mu _{n}^{G}\left( \varepsilon \right)
+\delta +\left( A+1\right) t\right\} \\
& \leq C_{1}N\left( \varepsilon \right) ^{2}\exp \left( -\frac{C_{2}\sqrt{n}\
\delta }{\left( N\left( \varepsilon \right) \right) ^{5/2}}\right) +2\exp
\left( -A_{1}\sqrt{n}t\right) +4\exp \left( -\frac{A_{5}t^{2}}{\varepsilon
^{2}}\right) .  
\end{split}
\end{equation*}%
\hfill $\square $

\begin{remark}
\label{rem:4} Here are the Polish spaces that allow us to apply the Berkes
and Philipp Lemma A1 as in the construction leading to (\ref{H9}). By
applying the entropy bound (\ref{eq:weak-bn1}) we can assume via the Dudley
entropy condition (\ref{DEC}) that the Gaussian process $\widetilde{\mathbb{F%
}}_{\left( \gamma _{n},T\right) }^{(n)}$ indexed by $\mathcal{F}_{n}$ in (%
\ref{H9}) is separable, bounded \ and $d_{P}$ uniformly continuous, where $%
d_{P}$ is defined as in (\ref{DP}). Moreover, since by using (\ref%
{eq:weak-bn1}), $\mathcal{F}_{n}$ is readily seen to be totally bounded, its
completion $\mathcal{F}_{n}^{c}$ is compact. (We complete $\mathcal{F}_{n}$
using the procedure described in Remark \ref{rem:c}.) Furthermore, the
process $\widetilde{\mathbb{F}}_{\left( \gamma _{n},T\right) }^{(n)}$ can be
readily extended to be a continuous Gaussian process on $\mathcal{F}_{n}^{c}$%
. Thus when applying the Berkes and Philipp lemma we can assume that $%
\widetilde{\mathbb{F}}_{\left( \gamma _{n},T\right) }^{(n)}$ is a Gaussian
process indexed by $\mathcal{F}_{n}^{c}$ taking values in the Polish space $%
S_{3}$ of bounded real valued functions defined on the compact set $\mathcal{%
F}_{n}^{c}$ continuous with respect to $d_{P}$. Therefore we can assume that 
$B_{1},\dots ,B_{n}$ i.i.d.~$B$, $Y_{1},\dots ,Y_{n}$ i.i.d.~$Y$ and $%
Z_{1},\dots ,Z_{n}$ i.i.d.~$Z$ take values in the Polish space $S_{1}\times
S_{2}$, where $S_{1}=\mathcal{C}\left( \left[ 0,T\right] \right) ^{n}\times 
\mathbb{R}^{N\left( \varepsilon \right) n}$ and $S_{2}=\mathbb{R}^{N\left(
\varepsilon \right) n}$, and $Z_{1},\dots, Z_{n}$ i.i.d.~$Z$ and $\widetilde{%
\mathbb{F}}_{\left( \gamma _{n},T\right) }^{(n)}$ take values in the Polish
space $S_{2}\times S_{3}$.
\end{remark}

The proof of Proposition \ref{propPrime} will be completed by refining
inequality (\ref{H9}). Recall the notation $\mathcal{F}\left( K,\gamma
\right) $, $\mathcal{G}\left( K,\gamma ,\varepsilon \right) $, (\ref{Ne})
and (\ref{notation}). We find that for any $0<\varepsilon ,u<e^{-1}$, with $%
K=M_{n}=\sqrt{c\log n}$, (\ref{bn1}) gives the bound, with $\nu_{0}$ and $%
H_{0}$ as in (\ref{nuk}), for some $c_{1}\geq 1$, 
\begin{equation}
\begin{split}
N(u)& =N_{[\;]}(u/2,\mathcal{F}\left( M_{n},\gamma \right) ,d_{P})
=N_{[\;]}(u/2,\mathcal{F}_{n},d_{P}) \\
& \leq c_{1}M_{n}^{1/H}u^{-\nu _{0}} \sqrt{\log u^{-1}}\gamma ^{-H_{0}}\left[
\log \left( M_{n}/(u\gamma )\right) \right] ^{1/H}
\end{split}
\label{H1}
\end{equation}%
and the weaker entropy bound (\ref{eq:weak-bn1}) combined with (\ref{NN})
implies that for some $c_{2}>0$ and any $u\in (0,1)$, $\gamma \in (0,1)$, 
\begin{equation}
\begin{split}
N_{[\;]}(u,\mathcal{G}\left( M_{n},\gamma ,\varepsilon \right) ,d_{P})&
=N_{[\;]}(u,\mathcal{G}_{n}\left( \varepsilon \right) ,d_{P}) \\
& \leq \left( N_{[\;]}(u/2,\mathcal{F}\left( M_{n},\gamma \right)
,d_{P})\right) ^{2} \\
& \leq c_{2}M_{n}^{4/H}u^{-3\nu _{0}}\gamma ^{-(2+4/H)}.
\end{split}
\label{H2}
\end{equation}

\noindent We shall make frequent use of the following elementary inequality.
For any $x\geq 1$ and any $\varepsilon \leq 1$ we have 
\begin{equation}
\int_{0}^{\varepsilon }\sqrt{x+\log u^{-1}}du\leq 2\varepsilon \sqrt{x+\log
\varepsilon ^{-1}}.  \label{ei}
\end{equation}%
Setting $\sigma =\varepsilon $ in (\ref{j}) and (\ref{a}) below we get using
(\ref{H2}) and (\ref{ei}) that for some $c_{3}>0$,
\begin{equation*}
J\left( \varepsilon ,\mathcal{G}_{n}\left( \varepsilon \right) \right)
=\int_{\left[ 0,\varepsilon \right] }\sqrt{1+\log N_{[\;]}(s,\mathcal{G}%
_{n}\left( \varepsilon \right) ,d_{P})}\,ds \leq c_{3}\,\varepsilon \,\sqrt{%
\log \left[ (\log n)/(\varepsilon \gamma )\right] }
\end{equation*}%
and for some $b_{0}>0$ 
\begin{equation*}
a\left( \varepsilon ,\mathcal{G}_{n}\left( \varepsilon \right) \right)
=\varepsilon \left( 1+\log N_{[\;]}(\varepsilon ,\mathcal{G}_{n}\left(
\varepsilon \right) ,d_{P})\right) ^{-1/2}\geq b_{0}\varepsilon \left( \log %
\left[ (\log n)/(\varepsilon \gamma )\right] \right) ^{-1/2}.
\end{equation*}%
For the $\mu _{n}^{S}\left( \varepsilon \right) $ in (\ref{H10}) we obtain
by inequality (\ref{mom3}) with measurable envelope $G=1$ 
\begin{equation*}
\mu _{n}^{S}\left( \varepsilon \right) \leq c_{3}\varepsilon A_{0}\sqrt{\log %
\left[ (\log n)/(\varepsilon \gamma )\right] }+A_{0}\sqrt{n}\,1{\left\{ 1>%
\sqrt{n}b_{0}\varepsilon /\sqrt{\log \left[ (\log n)/(\varepsilon \gamma )%
\right] }\right\} ,}
\end{equation*}%
which as long as $1>\varepsilon =\varepsilon _{n}>0$ and $1\geq \gamma
=\gamma _{n}>0$ satisfy 
\begin{equation}
\frac{\sqrt{n}\varepsilon _{n}}{\sqrt{\log \left[ (\log n)/(\varepsilon
_{n}\gamma _{n})\right] }}\rightarrow \infty ,\text{ as }n\rightarrow \infty
,  \label{zero}
\end{equation}%
implies that for all large enough $n$ for a suitable $A_{1}^{\prime }>0$ 
\begin{equation}
\mu _{n}^{S}\left( \varepsilon _{n}\right) \leq A_{1}^{\prime }\varepsilon
_{n}\sqrt{\log \left[ (\log n)/(\varepsilon _{n}\gamma _{n})\right] }.
\label{log}
\end{equation}

Recall the definition of $\mu _{n}^{G}\left( \varepsilon \right) $ in (\ref%
{H10}). We get via the Gaussian moment bound (\ref{gaussmom}),
(\ref{eq:weak-bn1}) and inequality (\ref{ei}) that for all $0<$ $\varepsilon _{n}<1/e$ 
and appropriate $A_{2}^{\prime }$ and $A_{3}^{\prime }$
\begin{equation*}
\mu _{n}^{G}\left( \varepsilon _{n}\right) \leq A_{2}^{\prime
}\int_{0}^{\varepsilon _{n}}\sqrt{\log \left[ (\log n)/(u\gamma _{n})\right] 
}\,du\leq A_{3}^{\prime }\varepsilon _{n}\sqrt{\log \left[ (\log
n)/(\varepsilon _{n}\gamma _{n})\right] }.
\end{equation*}%
Hence, as long as (\ref{zero}) is satisfied, for some $D>0$ we have for all
large enough $n$ 
\begin{equation}
A\mu _{n}^{S}\left( \varepsilon _{n}\right) +\mu _{n}^{G}\left( \varepsilon
_{n}\right) \leq D\varepsilon _{n}\sqrt{\log \left[ (\log n)/(\varepsilon
_{n}\gamma _{n})\right] }.  \label{eps}
\end{equation}%
In addition, by (\ref{H1}) we have the bound 
\begin{equation*}
N\left( \varepsilon _{n}\right) \leq c_{1}(\log n)^{1/(2H)}\varepsilon
_{n}^{-\nu _{0}}\sqrt{\log \varepsilon _{n}^{-1}}\gamma _{n}^{-H_{0}}\left(
\log \left[ (\log n)/(\varepsilon _{n}\gamma _{n})\right] \right) ^{1/H},
\end{equation*}%
and also the weaker bound (\ref{eq:weak-bn1}) gives for some $c_{4}>0$ 
\begin{equation*}
N\left( \varepsilon _{n}\right) \leq c_{4}(\log n)^{1/H}\varepsilon
_{n}^{-3(1+1/H)}\gamma _{n}^{-(1+2/H)}.
\end{equation*}%
Therefore, in view of (\ref{eps}) and (\ref{H9}), it is natural to define
for suitably large positive $\gamma _{1}^{\prime }$ and $\gamma _{2}^{\prime
}$, 
\begin{equation*}
\delta =\gamma _{1}^{\prime }\varepsilon _{n}\sqrt{\log \left[ (\log
n)/(\varepsilon _{n}\gamma _{n})\right] }\quad \text{and }\ t=\gamma
_{2}^{\prime }\varepsilon _{n}\sqrt{\log \left[ (\log n)/(\varepsilon
_{n}\gamma _{n})\right] }.
\end{equation*}%
We now have by (\ref{H9}), as long as (\ref{zero}) holds, that for all large
enough $n$ there is a probability space depending on $\gamma _{1}^{\prime
},\gamma _{2}^{\prime },$ $\gamma _{n}$ and $\varepsilon _{n}$ on which $%
\alpha _{n}$ and $\widetilde{\mathbb{F}}_{\left( \gamma _{n},T\right)
}^{(n)} $ sit such that 
\begin{equation*}
\begin{split}
& P\left\{ \left\Vert \alpha _{n}-\widetilde{\mathbb{F}}_{\left( \gamma
_{n},T\right) }^{(n)}\right\Vert _{\mathcal{F}_{n}}>\left( D+\gamma
_{1}^{\prime }+\left( 1+A\right) \gamma _{2}^{\prime }\right) \varepsilon
_{n}\sqrt{\log \left[ (\log n)/(\varepsilon _{n}\gamma _{n})\right] }\right\}
\\
& \leq \frac{C_{1}c_{4}^{2}(\log n)^{2/H}}{\gamma _{n}^{2+4/H}\varepsilon
_{n}^{3\nu _{0}}}\exp \left( -\frac{c_{5}\sqrt{n}\gamma _{1}^{\prime
}\varepsilon _{n}^{1+5\nu _{0}/2}\gamma _{n}^{5H_{0}/2}\sqrt{\log \left[
(\log n)/(\varepsilon _{n}\gamma _{n})\right] }}{(\log n)^{\frac{5}{4H}%
}(\log \varepsilon _{n}^{-1})^{5/4}\left( \log \left[ (\log n)/(\varepsilon
_{n}\gamma _{n})\right] \right) ^{\frac{5}{2H}}}\right) \\
& \phantom{\leq}+2\exp \left( -A_{1}\gamma _{2}^{\prime }\sqrt{n}\varepsilon
_{n}\sqrt{\log \left[ (\log n)/(\varepsilon _{n}\gamma _{n})\right] }\right)
+4\exp \left( -A_{5}(\gamma _{2}^{\prime })^{2}\log \left[ (\log
n)/(\varepsilon _{n}\gamma _{n})\right] \right) ,
\end{split}%
\end{equation*}%
for some $c_{5}>0$. Choose $\varepsilon _{n}$ such that 
\begin{equation*}
\sqrt{n}\varepsilon _{n}^{1+5\nu _{0}/2}\gamma _{n}^{5H_{0}/2}=(\log n)^{%
\frac{1}{2}+5\left( \frac{3}{4H}+\frac{1}{4}\right) }.
\end{equation*}%
Then by (\ref{eta}) 
\begin{equation*}
\frac{\log \varepsilon _{n}^{-1}}{\log n}\rightarrow \frac{1-5H_{0}\eta }{%
2+5\nu _{0}}>0,
\end{equation*}%
and 
\begin{equation*}
\frac{\log (\varepsilon _{n}\gamma _{n})^{-1}}{\log n}\rightarrow \frac{%
1-5H_{0}\eta }{2+5\nu _{0}}+\eta =:\zeta >0.
\end{equation*}%
An easy computation shows that the exponent of the first term satisfies with
a positive constant $\chi$ 
\begin{equation*}
\frac{\sqrt{n}\varepsilon _{n}^{1+5\nu _{0}/2}\gamma _{n}^{5H_{0}/2}\sqrt{%
\log \left[ (\log n)/(\varepsilon _{n}\gamma _{n})\right] }}{(\log n)^{\frac{%
5}{4H}}(\log \varepsilon _{n}^{-1})^{5/4}\left( \log \left[ (\log
n)/(\varepsilon _{n}\gamma _{n})\right] \right) ^{\frac{5}{2H}}}\sim \chi
\log n
\end{equation*}%
and 
\begin{equation*}
\varepsilon _{n}\sqrt{\log \left[ (\log n)/(\varepsilon _{n}\gamma _{n})%
\right] }\sim \sqrt{\zeta }\left( \log n\right) ^{\tau _{2}}\left(
n^{-1/2}\gamma _{n}^{-5H_{0}/2}\right) ^{1/(1+5\nu _{0}/2)},
\end{equation*}%
where $\tau _{2}$ is given in Proposition \ref{prop:1}. We readily obtain
from these last bounds that for every $\vartheta >1$ there exist $D>0$, $%
\gamma _{1}^{\prime }>0$ and $\gamma _{2}^{\prime }>0$ such that for all $%
n\geq 1$ large enough, $\alpha _{n}$ and $\widetilde{\mathbb{F}}_{\left(
\gamma _{n},T\right) }^{(n)}$ can be defined on the same probability space
so that 
\begin{equation*}
\begin{split}
P\bigg\{\left\Vert \alpha _{n}-\widetilde{\mathbb{F}}_{\left( \gamma
_{n},T\right) }^{(n)}\right\Vert _{\mathcal{F}_{n}}>& \left( D+\gamma
_{1}^{\prime }+\left( 1+A\right) \gamma _{2}^{\prime }\right) \\
& \times \sqrt{2\zeta }(\log n)^{\tau _{2}}\left( n^{-1/2}\gamma
_{n}^{-5H_{0}/2}\right) ^{1/(1+5\nu _{0}/2)}\bigg\}\leq n^{-\vartheta },
\end{split}%
\end{equation*}%
which in the special case when $\gamma _{n}=n^{-\eta }$, with $0\leq \eta <%
\frac{1}{5H_{0}}$, gives 
\begin{equation*}
P\left\{ \left\Vert \alpha _{n}-\mathit{\ }\widetilde{\mathbb{F}}_{\left(
\gamma _{n},T\right) }^{(n)}\right\Vert _{\mathcal{F}_{n}}>\left( D+\gamma
_{1}^{\prime }+\left( 1+A\right) \gamma _{2}^{\prime }\right) \sqrt{2\zeta }%
n^{-\tau _{1}}\left( \log n\right) ^{\tau _{2}}\right\} \leq n^{-\vartheta },
\end{equation*}%
where $\tau _{1}=\left( 1-5H_{0}\eta \right) /\left( 2+5\nu _{0}\right) >0$.
It is clear now that there exists a $\eta (\vartheta )>0$ such that (\ref{p}%
) and (\ref{pp}) hold. This completes the proof of Proposition \ref%
{propPrime}. \hfill $\square $ \smallskip

We are now ready to complete the proof of Proposition \ref{prop:1}. This
will be accomplished in two steps.\smallskip

\noindent \textbf{Step 1} We shall construct the needed version $\widetilde{%
\mathbb{G}}_{\left( \gamma _{n},T\right) }^{(n)}$ of the Gaussian process $%
\mathbb{G}_{\left( \gamma _{n},T\right) }$ as required in Proposition \ref%
{prop:1}. Set for a fixed $n$ for any $m\geq e$%
\begin{equation*}
\mathcal{F}_{m,n}=\left\{ h_{t,x}^{\left( m\right) }\left( g\right)
=1\left\{ g\left( t\right) \leq x,g\in \mathcal{C}_{m}\right\} :\left(
t,x\right) \in \mathcal{T}\left( \gamma _{n}\right) \right\}
\end{equation*}%
and let 
\begin{equation*}
\mathcal{F}_{\infty ,n}=\left\{ h_{t,x}^{\left( \infty \right) }\left(
g\right) =1\left\{ g\left( t\right) \leq x,g\in \mathcal{C}_{\infty
}\right\} :\left( t,x\right) \in \mathcal{T}\left( \gamma _{n}\right)
\right\} .
\end{equation*}%
(Note that $\mathcal{F}_{\infty ,n}=\mathcal{F}_{\left( \gamma_n ,T\right) }$%
.) Set 
\begin{equation*}
\mathcal{F}_{\infty }\left( \gamma _{n}\right) =\mathcal{F}_{\infty ,n}%
\mathcal{\cup }\cup _{m\geq e}\mathcal{F}_{m,n}.
\end{equation*}%
Let $\mathbb{H}_{\left( \gamma _{n},T\right) }$ be the mean zero Gaussian
process indexed by $\mathcal{F}_{\infty }\left( \gamma _{n}\right) $ such
that for $h_{s,x}^{(k)},h_{t,y}^{\left( m\right) }\in \mathcal{F}_{\infty
}\left( \gamma _{n}\right) $ with $e<k\leq m\leq \infty $ 
\begin{equation*}
\begin{split}
& \mathop{Cov}\left( \mathbb{H}_{\left( \gamma _{n},T\right) }\left(
h_{s,x}^{(k)}\right) ,\mathbb{H}_{\left( \gamma _{n},T\right) }\left(
h_{t,y}^{\left( m\right) }\right) \right) =P\left\{ B\left( s\right) \leq
x,B\left( t\right) \leq y,B\in \mathcal{C}_{k}\right\} \\
& \phantom{=} \ -P\left\{ B\left( s\right) \leq x,B\in \mathcal{C}%
_{k}\right\} P\left\{ B\left( t\right) \leq y,B\in \mathcal{C}_{m}\right\} .
\end{split}%
\end{equation*}%
In particular 
\begin{equation*}
\mathop{Cov}\left( \mathbb{H}_{\left( \gamma _{n},T\right) }\left(
h_{s,x}^{(\infty )}\right) ,\mathbb{H}_{\left( \gamma _{n},T\right) }\left(
h_{t,y}^{\left( \infty \right) }\right) \right) =\mathop{Cov}\left( \mathbb{G%
}_{\left( \gamma _{n},T\right) }\left( h_{s,x}\right), \mathbb{G}_{\left(
\gamma_{n},T\right) }\left( h_{t,y}\right) \right) .
\end{equation*}%
Thus $\mathbb{H}_{\left( \gamma _{n},T\right) }\left( h_{s,x}^{(\infty
)}\right)$ is a version of the Gaussian process $\mathbb{G}_{\left( \gamma
_{n},T\right) }\left( h_{s,x}\right)$. The process $\widetilde{\mathbb{G}}%
^{(n)}_{\left( \gamma_{n},T\right)}$ required in the statement of
Proposition 1 will be a version of the process $\mathbb{H}_{\left( \gamma
_{n},T\right) }\left( h_{s,x}^{(\infty )}\right)$. \smallskip

Notice that for $e<k\leq m\leq \infty $,%
\begin{equation}
\begin{split}
E\left( \mathbb{H}_{\left( \gamma _{n},T\right) }\left( h_{s,x}^{\left(
k\right) }\right) -\mathbb{H}_{\left( \gamma _{n},T\right) }\left(
h_{s,x}^{\left( m\right) }\right) \right) ^{2} & \leq P\left\{ B \left(
s\right) \leq x, B\in \mathcal{C}_{m}-\mathcal{C}_{k}\right\} \\
& \leq P\left\{ B\notin \mathcal{C}_{k}\right\} .
\end{split}
\label{H8a}
\end{equation}

In the following lemma using Dudley's entropy condition (\ref{Dd})
we show that $\mathbb{H}_{(\gamma_n,T)}$ has a 
continuous modification. To do so, we introduce further notation. 
For a set $\mathbb{T}$ equipped with
a semimetric $\rho$ let $N\left( \varepsilon ,\mathbb{T},\rho \right)$
denote the minimal number of $\rho$-balls of radius $\varepsilon$ needed
to cover $\mathbb{T}$. 

\begin{lemma} \label{lem:2} 
The Gaussian process $\mathbb{H}_{(\gamma _{n},T)}$ has a
bounded uniformly continuous modification 
$\widehat{\mathbb{H}}_{(\gamma_{n},T)}$.
\end{lemma}

\noindent \textit{Proof} 
Using the definition of $\mathcal{C}_k$ in (\ref{cn})
and the Landau--Shepp inequality (\ref{FLS1}) we obtain 
\begin{equation}  \label{eq:k}
P\{B\not\in \mathcal{C}_{k}\} = P \{ L > \sqrt{c \log k} \} 
\leq C e^{-D c\log k} = C k^{-Dc}.
\end{equation}
Let us fix $1 \geq \varepsilon > 0$ and choose $k = \lceil (4C /
\varepsilon^2)^{1/(cD)} \rceil$, where $\lceil \cdot \rceil$ stands for the
upper integer part. Then from (\ref{H8a}) and (\ref{eq:k}) follow for any $%
m\geq k$ (allowing $m=\infty $) that 
\begin{equation*}
d_{P}^{2}(h_{t,x}^{(m)},h_{t,x}^{(k)})\leq P\{B \not\in \mathcal{C}_{k}\}
\leq \varepsilon^2 /4.
\end{equation*}
For each $\ell \leq k$ choose a $d_{P}-\varepsilon /2$ grid $%
\{h_{t_{i},x_{i}}^{(\ell)}\}_{i=1}^{N_{\ell }}$ in $\mathcal{F}(\sqrt{%
c\log\ell}, \gamma_{n})=\mathcal{F}_{\ell,n}$. The entropy bound II
(\ref{eq:weak-bn1}) and the choice of $k$ shows that 
\begin{equation}
\begin{split}
N_{\ell }& \leq C(\log k)^{1/H}\varepsilon ^{-3(1+1/H)}\gamma _{n}^{-(1+2/H)}
\\
& \leq C^{\prime }\log \varepsilon ^{-1}\,\varepsilon ^{-3(1+1/H)}\gamma
_{n}^{-(1+2/H)}.
\end{split}
\label{eq:Nell-bound}
\end{equation}
Consider the finite set of functions 
\begin{equation*}
\mathcal{G}=\cup _{\ell \leq k}\{h_{t_{i},x_{i}}^{(\ell )}:i=1,2,\ldots
,N_{\ell }\}.
\end{equation*}
We claim that $\mathcal{G}$ is a $d_P-\varepsilon$ grid in $\mathcal{F}%
_\infty(\gamma_n)$. Indeed, let $h_{t,x}^{(m)} \in \mathcal{F}%
_\infty(\gamma_n)$ be arbitrary. If $m \leq k$ then there is an $%
h_{t_i,x_i}^{(m)} \in \mathcal{G}$ such that $d_P(h_{t,x}^{(m)},
h_{t_i,x_i}^{(m)} ) \leq \varepsilon/2$. For $m > k$ we have 
\begin{equation*}
d_P(h_{t,x}^{(m)}, h_{t_i,x_i}^{(k)}) \leq d_P(h_{t,x}^{(m)}, h_{t,x}^{(k)})
+ d_P(h_{t,x}^{(k)}, h_{t_i,x_i}^{(k)}) \leq \varepsilon/2 + \varepsilon/2=
\varepsilon,
\end{equation*}
where $h_{t_i,x_i}^{(k)}$ is chosen such that $d_P(h_{t,x}^{(k)},
h_{t_i,x_i}^{(k)} ) \leq \varepsilon/2$.

Thus $\mathcal{G}$ is indeed a $d_{P}-\varepsilon$ grid in
$\mathcal{F}_{\infty }(\gamma _{n})$, for which by (\ref{eq:Nell-bound}) 
\begin{equation}
N(\varepsilon, \mathcal{F}_\infty(\gamma_n), d_P) \leq 
|\mathcal{G}|=\sum_{\ell =1}^{k}N_{\ell }\leq C\varepsilon ^{-a}\gamma
_{n}^{-(1+2/H)},  \label{eq:G}
\end{equation}%
with $a=2/(cD)+6/H$, say.
Thus Dudley's condition (\ref{Dd}) is satisfied, and a bounded uniformly 
continuous modification $\widehat{\mathbb{H}}_{(\gamma _{n},T)}$ exists.
\hfill $\square$
\medskip

From now on to reduce notation we shall assume that
$\mathbb{H}_{(\gamma_{n},T)}$ is its bounded uniformly continuous modification. Consider 
the
class of functions $\mathcal{C}\left[ 0,T\right] \rightarrow \mathbb{R}^{2}$
indexed by $\left( t,x\right) \in \mathcal{T}\left( \gamma _{n}\right) $
given by 
\begin{equation*}
\mathcal{D}_{n}=\left\{ \left( h_{t,x}^{\left( n\right) },h_{t,x}^{\left(
\infty \right) }\right) :\left( t,x\right) \in \mathcal{T}\left( \gamma
_{n}\right) \right\}.
\end{equation*}%
Define the mean zero Gaussian process on $\mathcal{D}_{n}$ 
\begin{equation}
\mathbb{D}_{n}^{\left( n\right) }\left( h_{t,x}^{\left( n\right)
},h_{t,x}^{\left( \infty \right) }\right) =\mathbb{H}_{\left( \gamma
_{n},T\right) }\left( h_{t,x}^{\left( n\right) }\right) -\mathbb{H}_{\left(
\gamma _{n},T\right) }\left( h_{t,x}^{\left( \infty \right) }\right) .
\label{Z}
\end{equation}%
Introduce the semimetric on $\mathcal{D}_{n}$ 
\begin{equation*}
\rho _{P}^{\left( 1\right) }\left( \left( h_{s,x}^{\left( n\right)},
h_{s,x}^{\left( \infty \right) }\right),
\left( h_{t,y}^{\left( n\right)}, h_{t,y}^{\left( \infty \right) }\right) \right) 
=\sqrt{E\left(
\mathbb{D}_{n}^{\left( n\right) }
\left( h_{s,x}^{\left( n\right) },h_{s,x}^{\left(\infty \right) }\right) 
-\mathbb{D}_{n}^{\left( n\right) }\left(
h_{t,y}^{\left( n\right) },h_{t,y}^{\left( \infty \right) }\right) \right)^{2}}.
\end{equation*}%
Notice that%
\begin{equation*}
\begin{split}
\rho_{P}^{\left( 1\right) }\left( \left( h_{s,x}^{\left( n\right)
},h_{s,x}^{\left( \infty \right) }\right) ,\left( h_{t,y}^{\left( n\right)
},h_{t,y}^{\left( \infty \right) }\right) \right) 
& \leq \sqrt{2} \, d_{P}\left( h_{s,x}^{\left( n\right) },h_{t,y}^{\left( n\right) 
}\right) + \sqrt{2} \, d_{P}\left( h_{s,x}^{\left( \infty \right) },h_{t,y}^{\left(\infty 
\right) }\right) \\
&=:d_{P}^{\left( 1\right) }\left( \left(
h_{s,x}^{\left( n\right) },h_{s,x}^{\left( \infty \right) }\right) ,
\left(h_{t,y}^{\left( n\right) },h_{t,y}^{\left( \infty \right) }\right) \right) .
\end{split}
\end{equation*}
Thus $\rho _{P}^{\left( 1\right) }$ is bounded by the semimetric
$d_{P}^{\left( 1\right) }$. \smallskip 

With the view towards applying the Gaussian moment inequality (\ref{gm2}) let%
\begin{equation*}
\begin{split}
\mathop{diam}\left( \mathcal{D}_{n}\right)
& =\sup \left\{ \rho _{P}^{\left( 1\right)}
\left( \left( h_{s,x}^{\left( n\right) },h_{s,x}^{\left( \infty \right)
}\right) ,\left( h_{t,y}^{\left( n\right) },h_{t,y}^{\left( \infty \right)}
\right)\right) : \left( h_{s,x}^{\left( n\right) },
h_{s,x}^{\left( \infty \right) }\right) ,\left( h_{t,y}^{\left( n\right) },
h_{t,y}^{\left( \infty \right) }\right) \in \mathcal{D}_{n}\right\} \\
& =\sup \left\{ \rho _{P}^{\left( 1\right) }
\left( \left( h_{s,x}^{\left(n\right) }, h_{s,x}^{\left( \infty \right) }\right),
\left( h_{t,y}^{\left( n\right) }, h_{t,y}^{\left( \infty \right) }\right) \right)
:\left(s,x\right) ,(t,y)\in \mathcal{T}(\gamma _{n})\right\} 
\end{split}
\end{equation*}%
denote the diameter of the set $\mathcal{D}_{n}$. Observe that 
\begin{equation*}
\mathop{diam} \left( \mathcal{D}_{n}\right) \leq 
2\sup_{(s,x)\in \mathcal{T}(\gamma_{n})}
\sqrt{E\left( \mathbb{H}_{\left( \gamma _{n},T\right) }
\left( h_{s,x}^{\left( n\right) }\right) -\mathbb{H}_{\left( \gamma _{n},T\right)
}\left( h_{s,x}^{\left( \infty \right) }\right) \right) ^{2}},
\end{equation*}%
which by (\ref{H8a}) is 
\begin{equation*}
\leq 2\sqrt{P\left\{ B \notin \mathcal{C}_{n}\right\} }.
\end{equation*}%
This last bound, in turn, by the Landau--Shepp inequality (\ref{FLS1}) below
is 
\begin{equation}
=2\sqrt{P\left\{ L>\sqrt{c\log m}\right\} }\leq 2\sqrt{C}
\exp \left( -\frac{Dc\log n}{2}\right) .  \label{dd}
\end{equation}%
Thus for any $\Delta >1$ there exists a $c>0$ such that 
\begin{equation}
\mathop{diam}\left( \mathcal{D}_{n}\right) \leq n^{-\Delta }.  \label{ddd}
\end{equation}
Next notice that by the definition of $\mathcal{D}_n$ and
by (\ref{eq:G}) for some constant $c_{6}\geq 1$, 
\begin{equation}
\begin{split}
N(u,\mathcal{D}_{n},d_{P}^{\left( 1\right) })
\leq \left( N(u/(2\sqrt{2}), \mathcal{F}_{\infty }
\left( \gamma _{n}\right),d_{P} )\right)^{2} 
\leq c_{6}u^{-2a}\gamma _{n}^{-(2+4/H)}.
\end{split}
\label{nnu2}
\end{equation}%
Write%
\begin{equation*}
\left\Vert \mathbb{D}_{n}^{\left( n\right) }\right\Vert _{\mathcal{D}%
_{n}}=\sup \left\{ \left\vert \mathbb{D}_{n}^{\left( n\right) }\left(
h_{t,x}^{\left( n\right) },h_{t,x}^{\left( \infty \right) }\right)
\right\vert :\left( t,x\right) \in \mathcal{T}\left( \gamma _{n}\right)
\right\} .
\end{equation*}%
Combining (\ref{ddd}) and (\ref{nnu2}) with the Gaussian moment inequality
(\ref{gm2}) we have
\begin{equation*}
\begin{split}
\left\Vert \mathbb{D}_{n}^{\left( n\right) }\right\Vert _{\mathcal{D}_{n}}
& \leq 
E \left\vert \mathbb{H}_{\left( \gamma _{n},T\right) }
\left( h_{\gamma_{n},0}^{\left( n\right) }\right) -
\mathbb{H}_{\left( \gamma _{n},T\right) }
\left( h_{\gamma _{n},0}^{\left( \infty \right) }\right) \right\vert
+A_{4}\int_{0}^{n^{-\Delta }}
\sqrt{\log N(u,\mathcal{D}_{n},d_{P}^{\left(1\right) })}\,du \\
& \leq n^{-\Delta }+A_{4}\int_{0}^{n^{-\Delta }}\sqrt{\log c_{6}-2a\log
u-\left( 2+\frac{4}{H}\right) \log \gamma _{n}}\,du,
\end{split}
\end{equation*}%
which by using (\ref{eta}) and (\ref{ei}) gives for some $b>0$ , 
\begin{equation}
\left\Vert \mathbb{D}_{n}^{\left( n\right) }\right\Vert _{\mathcal{D}%
_{n}}\leq bn^{-\Delta }\sqrt{\log n}.  \label{bor1}
\end{equation}%
We have by using the Landau--Shepp inequality (\ref{FLS1}) that for a
sufficiently large $c>0$
\begin{equation}
\begin{split}
\sigma _{\mathcal{D}_{n}}^{2}\left( \mathbb{D}_{n}^{\left( n\right) }\right)
& =\sup \left\{ E\left( \mathbb{D}_{n}^{\left( n\right) }\left(
h_{t,x}^{\left( n\right) },h_{t,x}^{\left( \infty \right) }\right) \right)
^{2}:\left( t,x\right) \in \mathcal{T}\left( \gamma _{n}\right) \right\}  \\
& \leq P\left\{ B \notin \mathcal{C}_{n}\right\} \leq n^{-\Delta }.
\end{split}
\label{bor2}
\end{equation}%
Hence by Borell's inequality (\ref{bor}), for all $z>0$, 
\begin{equation*}
P\left\{ \left\vert \left\Vert \mathbb{D}_{n}^{\left( n\right) }\right\Vert
_{\mathcal{D}_{n}}-E\left\Vert \mathbb{D}_{n}^{\left( n\right) }\right\Vert
_{\mathcal{D}_{n}}\right\vert >z\right\} \leq 2\exp \left( -\frac{z^{2}}{%
2\sigma _{\mathcal{D}_{n}}^{2}\left( \mathbb{D}_{n}^{\left( n\right)
}\right) }\right) ,
\end{equation*}%
which on account of (\ref{bor1}) and (\ref{bor2}) gives for all $\theta >1$%
\begin{equation}
P\left\{ \left\Vert \mathbb{D}_{n}^{\left( n\right) }\right\Vert _{\mathcal{D%
}_{n}}>bn^{-\Delta }\sqrt{\log n}+2n^{-\Delta /2}\sqrt{\theta \log n}%
\right\} \leq n^{-\theta }.  \label{borbound}
\end{equation}

Returning to the construction of $\widetilde{\mathbb{G}}_{\left( \gamma
_{n},T\right) }^{(n)}$ in Proposition \ref{prop:1}, for each $n>e$, let $%
\mathbb{F}_{(\gamma _{n},T)}^{(n)}$ denote the restriction of $\mathbb{H}%
_{(\gamma _{n},T)}$ to $\mathcal{F}_{n}$ and $\mathbb{G}_{(\gamma _{n},T)}$
the restriction of $\mathbb{H}_{(\gamma _{n},T)}$ to $\mathcal{F}_{(\gamma
_{n},T)}$. Notice by (\ref{eq:Nell-bound}), $\mathcal{F}_{\infty }\left(
\gamma _{n}\right) $ is totally bounded in the $d_{P}$ semimetric, as are $%
\mathcal{F}_{n}$ and $\mathcal{F}_{\infty ,n}=\mathcal{F}_{\left( \gamma
,T\right) }$. By the discussion given in Remark \ref{rem:4} for $\mathbb{F}%
_{(\gamma _{n},T)}^{(n)}$ and $\mathcal{F}_{(\gamma _{n},T)}^{(n)}$, $%
\mathbb{G}_{\left( \gamma _{n},T\right) }$ can be extended to a continuous
function on the completion $\mathcal{F}_{\left( \gamma _{n},T\right) }^{c}$
of $\mathcal{F}_{\left( \gamma _{n},T\right) },$ which is compact. Therefore
we can argue that 
\begin{equation*}
\left( \mathbb{F}_{(\gamma _{n},T)}^{(n)}\text{,}\mathbb{G}_{\left( \gamma
_{n},T\right) },\left\Vert \mathbb{D}_{n}^{\left( n\right) }\right\Vert _{%
\mathcal{D}_{n}}\right) 
\end{equation*}%
takes values in the Polish space $S_{3}\times S_{4}$, where $S_{3}$ is as in
Remark \ref{rem:4} and $S_{4}=S_{3}^{\prime }\times \mathbb{R}$, with $%
S_{3}^{\prime }$ being the Banach space of bounded real valued functions
defined on the compact set $\mathcal{F}_{\left( \gamma _{n},T\right) }^{c}$
continuous with respect to $d_{P}$. Hence Lemma A1 of Berkes and Philipp
applies here and we can enlarge the probability space on which
inequality (\ref{p}) holds to include a version 
\begin{equation*}
\left( \widetilde{\mathbb{F}}_{(\gamma _{n},T)}^{(n)},
\widetilde{\mathbb{G}}_{\left( \gamma _{n},T\right) }^{\left( n\right) 
},\left\Vert \widetilde{\mathbb{D}}_{n}^{\left( n\right) }
\right\Vert _{\mathcal{D}_{n}}\right) 
\end{equation*}
of the process $\left( \mathbb{F}_{(\gamma _{n},T)}^{(n)},
\mathbb{G}_{\left( \gamma _{n},T\right) },
\left\Vert \mathbb{D}_{n}^{\left( n\right)} \right\Vert 
_{\mathcal{D}_{n}}\right)$ so that besides (\ref{p}), (\ref{borbound})
is also valid.\smallskip 

\noindent \textbf{Step 2}
We shall show that inequality (\ref{p1}) holds for 
$\left\Vert \alpha _{n}-
\widetilde{\mathbb{G}}_{\left( \gamma_{n},T\right)}^{\left(n\right)}
\right\Vert_{\mathcal{F}_{\left( \gamma _{n},T\right) }}$, which will complete 
the proof of Proposition \ref{prop:1}. Define for
$\left( h_{t,x}^{\left( n\right) },
h_{t,x}^{\left( \infty \right) }\right)\in \mathcal{D}_{n}$,
\begin{equation*}
\widetilde{\mathbb{D}}_{n}^{\left( n\right) }
\left( h_{t,x}^{\left( n\right)},
h_{t,x}^{\left( \infty \right) }\right) 
=\widetilde{\mathbb{F}}_{\left(\gamma _{n},T\right) }^{(n)}
\left( h_{t,x}^{(n)}\right) -
\widetilde{\mathbb{G}}_{\left( \gamma _{n},T\right)}^{\left( n\right) }
\left( h_{t,x}\right) .
\end{equation*}%
Clearly 
\begin{equation*}
\left\Vert \mathbb{D}_{n}^{\left( n\right) }\right\Vert _{\mathcal{D}_{n}}%
\overset{\mathrm{D}}{=}\left\Vert \widetilde{\mathbb{D}}_{n}^{\left(
n\right) }\right\Vert _{\mathcal{D}_{n}}=\sup \left\{ \left\vert \widetilde{%
\mathbb{F}}_{\left( \gamma _{n},T\right) }^{(n)}\left( h_{t,x}^{(n)}\right) -%
\widetilde{\mathbb{G}}_{\left( \gamma _{n},T\right) }^{\left( n\right)
}\left( h_{t,x}\right) \right\vert :\left( t,x\right) \in \mathcal{T}\left(
\gamma _{n}\right) \right\} .
\end{equation*}%
Notice that by (\ref{note}) 
\begin{equation*}
\begin{split}
\left\Vert \alpha _{n}-\widetilde{\mathbb{G}}_{\left( \gamma _{n},T\right)
}^{\left( n\right) }\right\Vert _{\mathcal{F}_{\left( \gamma _{n},T\right)
}}
& \leq \sup \left\{ \left\vert \alpha _{n}\left( h_{t,x}^{(n)}\right) -%
\widetilde{\mathbb{F}}_{\left( \gamma _{n},T\right) }^{(n)}\left(
h_{t,x}^{(n)}\right) \right\vert :\left( t,x\right) \in \mathcal{T}\left(
\gamma _{n}\right) \right\} \\
& \phantom{\leq} \ +\sum_{i=1}^{n}\frac{1\left\{ B_{i} \notin 
\mathcal{C}_{n}\right\} }
{\sqrt{n}}+\sqrt{n}P\left\{ B \notin \mathcal{C}_{n}\right\} +\left\Vert 
\widetilde{\mathbb{D}}_{n}^{\left( n\right) }\right\Vert _{\mathcal{D}_{n}}.
\end{split}
\end{equation*}%
Let 
\begin{equation*}
\delta _{n}\left( \Delta \right) =b n^{-\Delta }\sqrt{\log n}+
2n^{-\Delta /2} \sqrt{\theta \log n}.
\end{equation*}%
Recalling that $P\left\{ B \notin \mathcal{C}_{n}\right\} \leq n^{-\Delta }$, 
we get by (\ref{p}) and (\ref{borbound}) 
\begin{equation}
\begin{split}
& P\bigg\{\left\Vert \alpha _{n}-\widetilde{\mathbb{G}}_{\left( \gamma
_{n},T\right) }^{\left( n\right) }\right\Vert _{\mathcal{F}_{(\gamma
_{n},T)}}>\eta (\vartheta )\left( \log n\right) ^{\tau _{2}}\left(
n^{-1/2}\gamma _{n}^{-5H_{0}/2}\right) ^{2/(2+5\nu _{0})} \\
& \phantom{P \bigg\{ \left\Vert \alpha_{n} - \widetilde{\mathbb{G}}_{(
\gamma_{n},T ) } \right\Vert_{\mathcal{F}_{( \gamma_{n},T)}} > \ }+\sqrt{n}%
n^{-\Delta }+\delta _{n}(\Delta )\bigg\} \\
& \leq P\left\{ \left\Vert \alpha _{n}-\widetilde{\mathbb{F}}_{\left( \gamma
_{n},T\right) }^{(n)}\right\Vert _{\mathcal{F}_{n}}>\eta (\vartheta )\left(
\log n\right) ^{\tau _{2}}\left( n^{-1/2}\gamma _{n}^{-5H_{0}/2}\right)
^{2/(2+5\nu _{0})}\right\}  \\
& \phantom{\leq}+P\left\{ \sum_{i=1}^{n}1\left\{ B_{i}\notin \mathcal{C%
}_{n}\right\} >0\right\} +P\left\{ \left\Vert \mathbb{D}_{n}^{\left(
n\right) }\right\Vert _{\mathcal{D}_{n}}>\delta _{n}\left( \Delta \right)
\right\}  \\
& \leq n^{-\vartheta }+n^{1-\Delta }+n^{-\theta }.
\end{split}
\label{H4}
\end{equation}%
Noting that the $\Delta $ in the above inequalities can be made as large as
desired by choosing $c$ large enough, we see that for every $\lambda >1$,
for sufficiently large $c>0$ (that is large $\Delta >0$), $\vartheta >0$, $%
\theta >0$ and all large $n$%
\begin{equation}
n^{-\theta }+n^{1-\Delta }+n^{-\vartheta }<n^{-\lambda },  \label{omega}
\end{equation}%
and for any choice of $\vartheta >0$, $\theta >0$ and large enough $c>0$
(large $\Delta >0$), for all large $n$%
\begin{equation*}
\eta (\vartheta )\left( n^{-1/2}\gamma _{n}^{-5H_{0}/2}\right) ^{2/(2+5\nu
_{0})}\geq \sqrt{n}n^{-\Delta }+\delta _{n}\left( \Delta \right) .
\end{equation*}%
Thus there is a $\rho \left( \lambda \right) >0,$ and $c>0$, $\vartheta >0$
and $\theta >0$ such that for all large enough $n$,%
\begin{equation*}
\rho (\lambda )(\log n)^{\tau _{2}}\left( n^{-1/2}\gamma
_{n}^{-5H_{0}/2}\right) ^{2/(2+5\nu _{0})}
\end{equation*}%
\begin{equation*}
>\eta (\vartheta )\left( \log n\right) ^{\tau _{2}}\left( n^{-1/2}\gamma
_{n}^{-5H_{0}/2}\right) ^{2/(2+5\nu _{0})}+\sqrt{n}n^{-\Delta }+\delta
_{n}\left( \Delta \right) 
\end{equation*}%
and (\ref{omega}) holds, which by (\ref{H4}) implies that 
\begin{equation*}
P\left\{ \left\Vert \alpha _{n}-\widetilde{\mathbb{G}}_{\left( \gamma
_{n},T\right) }^{\left( n\right) }\right\Vert _{\mathcal{F}_{\left( \gamma
_{n},T\right) }}>\rho (\lambda )\left( \log n\right) ^{\tau _{2}}\left(
n^{-1/2}\gamma _{n}^{-5H_{0}/2}\right) ^{2/(2+5\nu _{0})}\right\}
<n^{-\lambda },
\end{equation*}%
that is, for all such large $n$ there exists a suitable probability space
such that (\ref{p1}) holds. This completes the proof of Proposition \ref{prop:1}.

\subsection{Proof of Proposition \protect\ref{prop:2}}

Put $\gamma _{n}=n^{-\eta }$, with $\eta =(5H_{0}+\kappa (2+5\nu _{0}))^{-1}$.
Note that for this choice of $\eta$ 
\begin{equation}
\tau_1\left( \eta \right) = \tau _{1}^{\prime }=\tau_1^{\prime}\left( \kappa
\right) =\kappa /(5H_{0}+\kappa (2+5\nu _{0}))=\kappa \eta .  \label{knu}
\end{equation}%
Applying Proposition \ref{prop:1}, for every $\lambda ^{\prime }>\lambda >1$
there exists a $\rho (\lambda ^{\prime })>0$ such that for each integer $n$
large enough one can construct on the same probability space random vectors $%
B_{1},\ldots ,B_{n}$ i.i.d.~$B$ and a probabilistically equivalent version $%
\widetilde{\mathbb{G}}_{\left( \gamma_{n},T\right) }^{(n)}$ of $\mathbb{G}%
_{\left( \gamma _{n},T\right) }$ such that, 
\begin{equation*}
P\left\{ \left\Vert \alpha _{n}-\widetilde{\mathbb{G}}_{\left( \gamma
_{n},T\right) }^{(n)}\right\Vert _{\mathcal{F}_{\left( \gamma _{n},T\right)
}}>\rho \left( \lambda ^{\prime }\right) n^{-\tau _{1}}\left( \log n\right)
^{\tau _{2}}\right\} \leq n^{-\lambda ^{\prime }},
\end{equation*}%
with $\tau _{1}=\tau _{1}\left( \eta \right) =\left( 1-5H_{0}\eta \right)
/\left( 2+5\nu _{0}\right) $, which, since $T^{\kappa }/t^{\kappa }\geq 1$
for $t\in \left[ \gamma _{n},T\right] $, implies that 
\begin{equation}
P\left\{ \sup_{\left( t,x\right) \in \left[ \gamma _{n},T\right] \times 
\mathbb{R}}t^{\kappa }\left\vert \alpha _{n}\left( h_{t,x}\right) -%
\widetilde{\mathbb{G}}_{\left( \gamma _{n},T\right) }^{(n)}\left(
h_{t,x}\right) \right\vert >T^{\kappa }\rho \left( \lambda ^{\prime }\right)
n^{-\tau _{1}}\left( \log n\right) ^{\tau _{2}}\right\} \leq n^{-\lambda
^{\prime }}.  \label{p11}
\end{equation}%
\smallskip

Using Lemma A1 of Berkes and Philipp, we can enlarge the probability on
which (\ref{p11}) holds to include a Gaussian process $\mathbb{G}_{\left(
0,T\right) }$ indexed by $\mathcal{G}\left( \kappa\right) $, so that $%
\mathbb{G}_{\left( 0,T\right) }$ and $\widetilde{\mathbb{G}}_{\left(
\gamma_{n},T\right) }^{(n)}$ agree on 
\begin{equation*}
\left\{ t^{\kappa}h_{t,x}:\left( t,x\right) \in\left[ \gamma_{n},T\right]
\times\mathbb{R}\right\} .
\end{equation*}
(The validity of the application of the Berkes and Philipp lemma can be
argued as in Remark \ref{rem:4}.) Further we have, using inequality (\ref%
{ginq}) below with $\delta=\kappa,$ that for a suitable $d_{1}>0$ for all
large $n$ 
\begin{equation}
P\left\{ \sup\left\{ t^{\kappa}\left\vert G\left( t,x\right) \right\vert
:\left( t,x\right) \in\left[ 0,n^{-\eta}\right] \times\mathbb{R}\right\}
>d_{1}n^{-\eta\kappa}\sqrt{\log n}\right\} \leq n^{-\lambda^{\prime}},
\label{ineqa}
\end{equation}
where $G\left( t,x\right) =\mathbb{G}_{\left( 0,T\right) }\left(
h_{t,x}\right) $ for $\left( t,x\right) \in\left[ 0,n^{-\eta}\right] \times%
\mathbb{R}$.\smallskip

Next by using inequality (\ref{videlta}) below with $\delta=\kappa$ we get
that for a suitable $d_{2}>0$ for all large $n$ 
\begin{equation}
P\left\{ \sup\left\{ t^{\kappa}\left\vert \alpha_{n}\left( h_{t,x}\right)
\right\vert :\left( t,x\right) \in\left[ 0,n^{-\eta}\right] \times \mathbb{R}%
\right\} >d_{2}n^{-\eta\kappa}\sqrt{\log n}\right\} \leq
n^{-\lambda^{\prime}}.  \label{ineqb}
\end{equation}
Recall the notation after (\ref{not}). Combining inequalities (\ref{p11}), (\ref%
{ineqa}) and (\ref{ineqb}), and noting that $\tau_{2}>1/2$, we get for all
large enough $n$ 
\begin{equation*}
P\left\{ \left\Vert \alpha_{n}-\widetilde{\mathbb{G}}_{\left( 0,T\right)
}\right\Vert _{\mathcal{G}(\kappa)}>\left( d_{1}+d_{2}+T^{\kappa}\rho
(\lambda^{\prime})\right) n^{-\min\{\tau_{1},\eta\kappa\}}\left( \log
n\right) ^{\tau_{2}}\right\} \leq3n^{-\lambda^{\prime}}.
\end{equation*}
It is clear now that the optimal choice for $\eta$ satisfies $%
\tau_{1}(\eta)=\kappa\eta$, which by (\ref{knu}) our chosen value fulfills.
Thus by choosing $\lambda^{\prime}$ so that $3n^{-\lambda^{\prime}}<n^{-%
\lambda}$, setting $\rho^{\prime}\left( \lambda\right)
=d_{1}+d_{2}+T^{\kappa}\rho\left( \lambda^{\prime}\right) $, we conclude
that (\ref{del}) holds. \hfill$\square$

\begin{remark}
\label{rem:6} Here we discuss the continuity of the Gaussian process $%
\mathbb{G}_{\left( 0,T\right) }$ indexed by $\mathcal{G}\left( \kappa
\right) $. A straightforward argument based on Inequality 1 in the Appendix
shows that, w.p.~$1$, for all $\varepsilon >0$ there exists a $0<\gamma <1$
such that 
\begin{equation}
\sup_{(t,x)\in \lbrack 0,\gamma ]\times \mathbb{R}}t^{\kappa }\left\vert
G\left( t,x\right) \right\vert <\varepsilon .  \label{gdelta}
\end{equation}%
Moreover, as pointed out above, for any $0<\gamma <1$, $G\left( t,x\right)$
is almost surely bounded and uniformly continuous on $[\gamma ,T]\times 
\mathbb{R}$, when equipped with the semimetric (\ref{metric}), which implies
the same for $t^{\kappa }G\left( t,x\right) $, which when combined with (\ref%
{gdelta}), readily implies that $t^{\kappa }G\left( t,x\right) $ is almost
surely bounded and uniformly continuous on $\mathcal{T}\left( 0\right)
=[0,T]\times \mathbb{R}$ with respect to the semimetric $\rho_{\kappa }$ defined in 
(\ref{mk}). Also by applying Proposition 1 on page 26 of Lifshits \cite{Lifshits} we can 
assume that the Gaussian process $t^{\kappa }G\left(t,x\right) $ is separable. Thus there 
exists a version of $t^{\kappa}G\left( t,x\right)$ that is bounded and uniformly 
continuous on $\mathcal{T}\left( 0\right)$.
\end{remark}

\section{Appendix}
\label{s6}

\subsection{A Gaussian coupling inequality}

\noindent Einmahl and Mason \cite{EM97} pointed out in their Fact 2.2 that
the Strassen--Dudley theorem (see Theorem 11.6.2 in Dudley \cite{DU89}) in
combination with a special case of Theorem 1.1 of Zaitsev \cite{87a} (also
see the discussion after its statement) yields the following Gaussian
coupling. Here $\left\vert \cdot\right\vert_{N}$, $N\geq1$, denotes the
usual Euclidean norm on $\mathbb{R}^{N}$. \smallskip

\noindent \textbf{Coupling inequality.} \textit{Let $Y_{1},\ldots ,Y_{n}$ be
independent mean zero random vectors in $\mathbb{R}^{N}$, $N\geq 1$, such
that for some $b>0$, 
\begin{equation*}
\left\vert Y_{i}\right\vert _{N}\leq b,\ i=1,\dots ,n.
\end{equation*}%
If $(\Omega, \mathcal{T}, P)$ is rich enough then for each $\delta >0$, one
can define independent normally distributed mean zero random vectors $%
Z_{1},\ldots ,Z_{n}$ with $Z_{i}$ and $Y_{i}$ having the same covariance
matrix for $i=1,\ldots ,n$, such that for universal constants $C_{1}>0$ and $%
C_{2}>0$, 
\begin{equation}
P\left\{ \left\vert \sum_{i=1}^{n}\left( Y_{i}-Z_{i}\right) \right\vert
_{N}>\delta \right\} \leq C_{1}N^{2}\exp \left( -\frac{C_{2}\delta }{N^{2}b}%
\right) .  \label{coup}
\end{equation}%
}

\begin{remark}
\label{rem:5} Actually Einmahl and Mason did not specify the $N^{2}$ in (\ref%
{coup}) and they applied a less precise result given Theorem 1.1 in \cite%
{87b} with $N^{2}$ replaced by $N^{5/2}$, however their argument is equally
valid when based upon Theorem 1.1 in \cite{87a}. Zaitsev \cite{87a} remarks
that the assumptions of Theorem 1.1 of \cite{87b} imply those of Theorem 1.1
of \cite{87a}. See, in particular, the paragraph right above Remark 1.1 in 
\cite{87a}. Also see equation (18) in \cite{ZZ}.
\end{remark}

\subsection{Pointwise measurable classes}

\noindent \textbf{Definition.} A class $\mathcal{G}$ of measurable
real-valued functions defined on a measurable space $\left( S,\mathcal{S}%
\right) $ is \textit{pointwise measurable} if there exists a countable
subclass $\mathcal{G}_{\infty }$\ of $\mathcal{G}$ such that we can find for
any function $f\in \mathcal{G}$ a sequence of functions $\{f_{m}\}$ in $%
\mathcal{G}_{\infty }$\ for which $\lim_{m\rightarrow \infty }f_{m}(x)=f(x)$%
\ for all $x\in S$. For more about pointwise measurability see pages 109-110
and Example 2.3.4 of van der Vaart and Wellner \cite{VaartWellner}, as well
as Section 8.2 of Kosorok \cite{koro}.\smallskip

We shall show here that the classes of functions $\mathcal{F}\left( K,\gamma
\right) $, $K\geq 1$, of the form (\ref{FK}), where $0\leq \gamma
<1<T<\infty $ are pointwise measurable. Let $\mathbb{Q}$ denote the set of
rational numbers. For any $K\geq 1$ consider the countable class $\mathcal{F}%
_{\infty ,K}$ of functions of $g\in \mathcal{C}\left[ 0,T\right] \rightarrow
\left\{ 0,1\right\} $ indexed by $u,v\in \left[ \gamma ,T\right] \cap 
\mathbb{Q\cup }\left\{ \gamma ,T\right\}, y\in \mathbb{Q}$ defined by 
\begin{equation*}
1\left\{ g\left( v\right) -Kf_{H}(\left\vert v-u\right\vert )\leq y,\text{ }%
g\in \mathcal{C}\left( K\right)\right\},
\end{equation*}%
where $\mathcal{C}\left( K\right) $ is as in (\ref{cck}). Clearly for each $%
\left( t,x\right) \in \mathcal{T}(\gamma )=\left[ \gamma ,T\right] \times 
\mathbb{R}\ $we can choose sequences $s_{m}$ and $t_{m}\in \left[ \gamma ,T%
\right] \cap \mathbb{Q\cup }\left\{ \gamma ,T\right\} $ such that $%
t_{m}\searrow t$ and $s_{m}\nearrow t$. Also we can select a sequence $%
y_{m}\in \mathbb{Q}\searrow x.$ We see that each 
\begin{equation*}
1\left\{ g\left( t_{m}\right) -Kf_{H}(\left\vert t_{m}-s_{m}\right\vert
)\leq y_{m},g\in \mathcal{C}\left( K\right) \right\} \in \mathcal{F}_{\infty
,K}\text{.}
\end{equation*}%
Moreover, if $g\in \mathcal{C}\left( K\right) $, then $g\left( t_{m}\right)
-Kf_{H}(\left\vert t_{m}-s_{m}\right\vert )\leq g\left( t\right) $ and $%
g\left( t_{m}\right) -Kf_{H}(\left\vert t_{m}-s_{m}\right\vert )\rightarrow
g\left( t\right) $. Thus if $g\left( t\right) \leq x$ and $g\in \mathcal{C}%
\left( K\right) $ then 
\begin{equation*}
1\left\{ g\left( t_{m}\right) -Kf_{H}(\left\vert t_{m}-s_{m}\right\vert
)\leq y_{m},g\in \mathcal{C}\left( K\right) \right\} =1\rightarrow
1=h_{t,x}^{\left( K\right) }\left( g\right).
\end{equation*}%
$\ $Whereas if $g\left( t\right) >x$ then for some $\delta >0$, $g\left(
t\right) >x+\delta $ and all large enough $m$, 
\begin{equation*}
g\left( t_{m}\right) -Kf_{H}(\left\vert t_{m}-s_{m}\right\vert )>x+\delta /2%
\text{ and }x+\delta /4>y_{m}.
\end{equation*}%
This says that eventually $g\left( t_{m}\right) -Kf_{H}(\left\vert
t_{m}-s_{m}\right\vert )>y_{m}$ and thus 
\begin{equation*}
1\left\{ g\left( t_{m}\right) -Kf_{H}(\left\vert t_{m}-s_{m}\right\vert
)\leq y_{m},g\in \mathcal{C}\left( K\right) \right\} =0=h_{t,x}^{\left(
K\right) }\left( g\right) .
\end{equation*}%
Hence $\mathcal{F}\left( K,\gamma \right) $ is pointwise measurable with
countable subclass $\mathcal{F}_{\infty ,K}$.\medskip

For any $\kappa>0$ and $K\geq1$ let $\mathcal{G}\left( \kappa,K\right) $
denote the class of functions $g\in\mathcal{C}\left[ 0,T\right] \rightarrow%
\left[ 0,T^{\kappa}\right] $ indexed by $\left( t,x\right) \in\mathcal{T}%
\left( 0\right) =\left[ 0,T\right] \times\mathbb{R}$ defined by 
\begin{equation}
t^{\kappa}h_{t,x}^{\left( K\right) }\left( g\right) =t^{\kappa}1\left\{
g\left( t\right) \leq x,g\in\mathcal{C}\left( K\right) \right\} .  \label{gk}
\end{equation}
Clearly by a slight modification of the above argument $\mathcal{G}\left(
\kappa,K\right) $ is pointwise measurable.

\subsection{Inequalities for empirical processes}

In this subsection $\mathcal{G}$ is a pointwise measurable class of
measurable real-valued functions defined on a measurable space $\left( S,%
\mathcal{S}\right) $. For any $0<\sigma <1$, set 
\begin{equation}
J\left( \sigma ,\mathcal{G}\right) =\int_{\left[ 0,\sigma \right] }\sqrt{%
1+\log N_{[\;]}(s,\mathcal{G},d_{P})}\,ds  \label{j}
\end{equation}%
and 
\begin{equation}
a\left( \sigma ,\mathcal{G}\right) =\sigma \left[ 1+\log N_{[\;]}(\sigma ,%
\mathcal{G},d_{P})\right] ^{-1/2}.  \label{a}
\end{equation}%
Lemma 19.34 in van der Vaart \cite{Vaart} gives the following moment bound.
(Note the needed \textquotedblleft\ $+1$\textquotedblright\ in the
definition of $J(\sigma ,\mathcal{G})$ and $a\left( \sigma ,\mathcal{G}%
\right) $.)\medskip

\noindent \textbf{Moment inequality.} \textit{Let $\xi ,\xi
_{1},\ldots,\xi_{n}$ be i.i.d.~and assume that $\mathcal{G}$ has a
measurable envelope function $G$ and $E\left( g^{2}\left( \xi \right)
\right) <\sigma^{2}<1$ for every $g\in \mathcal{G}$. We have, for a
universal constant $A_{0}^{\prime}$, 
\begin{equation}
\begin{split}
E \left\Vert \frac{1}{\sqrt{n}}\sum_{i=1}^{n}\left( g(\xi _{i})-Eg(\xi
_{i})\right) \right\Vert _{\mathcal{G}} \leq A_{0}^{\prime }\left[ J\left(
\sigma ,\mathcal{G}\right) +\sqrt{n}\ E\left( G\left( \xi \right) 1\left\{
G\left( \xi \right) >\sqrt{n}a(\sigma, \mathcal{G})\right\} \right) \right] .
\end{split}
\label{H5}
\end{equation}%
}

Let $\epsilon $ be a Rademacher variable, i.e.~$P\{\epsilon =1\}=P\{\epsilon
=-1\}=1/2$, and consider independent Rademacher variables $\epsilon
_{1},\ldots ,\epsilon _{n}$ independent of $\xi _{1}$, $\ldots ,\xi _{n}$.
From a special case of a well-known symmetrization lemma, we have for any
class of functions $\mathcal{G}$ in $L_{1}\left( P\right) $ 
\begin{equation*}
\begin{split}
\frac{1}{2}E\left\Vert \sum_{i=1}^{n}\epsilon _{i}\left( g(\xi
_{i})-Eg\left( \xi \right) \right) \right\Vert _{\mathcal{G}} \leq
E\left\Vert \sum_{i=1}^{n}\left( g(\xi _{i})-Eg\left( \xi \right) \right)
\right\Vert _{\mathcal{G}} \leq 2E\left\Vert \sum_{i=1}^{n}\epsilon
_{i}g(\xi_{i}) \right\Vert_{\mathcal{G}}.
\end{split}%
\end{equation*}%
(See Lemma 6.3 of Ledoux and Talagrand \cite{LedeouxTalagrand3}.) In
particular we get
\begin{equation}
\begin{split}
E\left\Vert \sum_{i=1}^{n}\epsilon _{i}g(\xi _{i})\right\Vert _{\mathcal{G}%
}& \leq E\left\Vert \sum_{i=1}^{n}\epsilon _{i}\left( g(\xi _{i})-Eg\left(
\xi \right) \right) \right\Vert _{\mathcal{G}}+E\left\vert \sum_{i=1}^{n}%
\mathbb{\epsilon }_{i}\right\vert \left\Vert E g\left( \xi \right)
\right\Vert _{\mathcal{G}}  \\
& \leq 2E\left\Vert \sum_{i=1}^{n}\left( g(\xi _{i})-Eg\left( \xi \right)
\right) \right\Vert _{\mathcal{G}}+\sigma \sqrt{n}.  \label{H5A}
\end{split}
\end{equation}
Thus we readily get from (\ref{H5}) with $A_{0}=2A_{0}^{\prime }+1$ and
noting that the integrand of $J\left( \sigma ,\mathcal{G}\right) $ is
greater than or equal to $1$, 
\begin{equation}
E\left\Vert \frac{1}{\sqrt{n}}\sum_{i=1}^{n}\epsilon _{i}g(X_{i})\right\Vert
_{\mathcal{G}}\leq A_{0}\left[ J\left( \sigma ,\mathcal{G}\right) +\sqrt{n}\
E\left( G\left( \xi \right) 1\left\{ G\left( \xi \right) >\sqrt{n}\ a(\sigma
,\mathcal{G})\right\} \right) \right] .  \label{mom3}
\end{equation}

We shall be using the moment bound (\ref{mom3}) in conjunction with the
following exponential inequality due to Talagrand \cite{Talagrand}. This
maximal version is pointed out by Einmahl and Mason \cite[Inequality A.1 on
p.31]{EinmahlMason}.\medskip

\noindent \textbf{Talagrand inequality.} \textit{Let $\mathcal{G}$ be a
pointwise measurable class of measurable real-valued functions defined on a
measurable space $(S,\mathcal{S})$ satisfying $||g||_{\infty }\leq M,\ g\in 
\mathcal{G}$, for some $0<M<\infty $. Let $X,X_{n}$, $n\geq 1$, be a
sequence of i.i.d.~random variables defined on a probability space $\left(
\Omega ,\mathcal{A},P\right) $ and taking values in $S$, then for all $t>0$
we have for suitable finite constants $A,A_{1}>0$, 
\begin{equation}
\begin{split}
P\bigg\{ \max_{1\leq m\leq n} \|\sqrt{m}\alpha _{m}\|_{\mathcal{G}}\geq
A \bigg( E \, \Big\Vert \sum_{i=1}^{n}\epsilon _{i}g(X_{i}) 
\Big\Vert_{\mathcal{G}}+t \bigg) \bigg\} \leq 2\exp \bigg( \!\!-\frac{A_{1}t^{2}}{n\sigma
_{\mathcal{G}}^{2}} \bigg)\! +2\exp \bigg(\!\! -\frac{A_{1}t}{M} \bigg) ,
\end{split}
\label{tal}
\end{equation}%
where $\sigma _{\mathcal{G}}^{2}=\sup_{g\in \mathcal{G}}\mathop{Var}(g(X))$. 
}

\subsection{Inequalities for Gaussian processes}

Let $\mathbb{Z}$ be a separable mean zero Gaussian process on a probability
space $(\Omega ,\mathcal{A},P)$ indexed by a set $\mathbb{T}$, equipped with
a semimetric 
\begin{equation}
\rho \left( s,t\right) =\sqrt{E \left( \mathbb{Z}\left( t\right) -\mathbb{Z}%
\left( s\right) \right) ^{2}}.  \label{rho1}
\end{equation}%
For each $\varepsilon >0$ let $N\left( \varepsilon ,\mathbb{T},\rho \right) $
denote the minimal number of $\rho $-balls of radius $\varepsilon $ needed
to cover $\mathbb{T}.$ Write $\left\Vert \mathbb{Z}\right\Vert _{\mathbb{T}%
}=\sup_{t\in \mathbb{T}}\left\vert \mathbb{Z}_{t}\right\vert $ and $\sigma _{%
\mathbb{T}}^{2}\left( \mathbb{Z}\right) =\sup_{t\in \mathbb{T}}E\left( 
\mathbb{Z}_{t}^{2}\right) $.\smallskip

According to Dudley \cite{Dudley}, the entropy condition 
\begin{equation}
\int_{\left[ 0,1\right] }\sqrt{\log N\left( \varepsilon ,\mathbb{T},\rho
\right) }\,d\varepsilon <\infty  \label{DEC}
\end{equation}%
ensures the existence of a separable, bounded, $\rho $-uniformly continuous
modification of $\mathbb{Z}$. The following moment bound is a version of
Corollary 2.2.8 in van der Vaart and Wellner \cite{VaartWellner}. (Also see
their Problem 2.2.14.)\smallskip

\noindent \textbf{Gaussian moment inequality.} \textit{For some universal
constant $A_{4}>0$ and all $\sigma >0$ we have 
\begin{equation}
E\left( \sup_{\rho \left( s,t\right) <\sigma }\left\vert \mathbb{Z}_{t}-%
\mathbb{Z}_{s}\right\vert \right) \leq A_{4}\int_{\left[ 0,\sigma \right] }%
\sqrt{\log N\left( \varepsilon ,\mathbb{T},\rho \right) }\,d\varepsilon
\label{gaussmom}
\end{equation}%
and for any $t_{0}\in \mathbb{T}$, 
\begin{equation}
E\left( \left\Vert \mathbb{Z}\right\Vert _{\mathbb{T}}\right) \leq
E\left\vert \mathbb{Z}_{t_{0}}\right\vert +A_{4}\int_{\left[ 0,\mathbb{D}%
\right] }\sqrt{\log N\left( \varepsilon ,\mathbb{T},\rho \right) }%
\,d\varepsilon ,  \label{gm2}
\end{equation}
with 
\begin{equation}
\mathbb{D}=\sup_{s,t\in \mathbb{T}}\rho \left( s,t\right)  \label{dia}
\end{equation}
denoting the diameter of $\mathbb{T}$. } \smallskip

\noindent Notice that if $d$ is a semimetric on $\mathbb{T}$ such that for
all $s,t\in T$, $d\left( s,t\right) \geq \rho \left( s,t\right) $, then 
\begin{equation*}
\sup_{\left\{ s:\rho \left( s,t\right) <\sigma \right\} }\left\vert \mathbb{Z%
}_{t}-\mathbb{Z}_{s}\right\vert \geq \sup_{\left\{ s:\text{ }d\left(
s,t\right) <\sigma \right\} }\left\vert \mathbb{Z}_{t}-\mathbb{Z}%
_{s}\right\vert
\end{equation*}%
and $N\left( \varepsilon ,\mathbb{T},d\right) \geq $ $N\left( \varepsilon ,%
\mathbb{T},\rho \right) $. Thus 
\begin{equation}
\int_{\left[ 0,1\right] }\sqrt{\log N\left( \varepsilon ,\mathbb{T},d\right) 
}\,d\varepsilon <\infty  \label{Dd}
\end{equation}%
implies by the Dudley result the existence of a separable, bounded, $d$%
-uniformly continuous modification of $\mathbb{Z}$. (Here note that $\rho $%
-uniformly continuous implies $d$-uniformly continuous.) Moreover the moment
inequalities in (\ref{gaussmom}) and (\ref{gm2}) hold when $\rho $ is
replaced by $d$ and in the definition of $\mathbb{D}.$\smallskip

In particular, these inequalities hold when $\mathbb{Z}=\mathbb{G}%
_{(\gamma,T)}$, the Gaussian process defined at the end of Subsection 2.1,
where $\mathbb{T}=\mathcal{F}_{(\gamma,T)}$ and $d=d_{P}$ is as defined in (%
\ref{DP}), and $\mathbb{D}=\sup \left\{d_P \left(f,g\right): f,g\in \mathcal{%
F}_{(\gamma,T)}\right\}$ is the diameter $\mathbb{D}$ of $\mathbb{T}=%
\mathcal{F}_{(\gamma,T)}$. \smallskip

The following large deviation probability estimate for $\left\Vert \mathbb{Z}%
\right\Vert _{\mathbb{T}}$ is due to Borell \cite{BorellChrister}. (Also see
Proposition A.2.1 in \cite{VaartWellner}.) Let $M\left( X\right) $ denote
the a median of $\left\Vert \mathbb{Z}\right\Vert _{\mathbb{T}}$, i.e.~$%
P\left\{ \left\Vert \mathbb{Z}\right\Vert _{\mathbb{T}}\geq M\left( X\right)
\right\} \geq1/2$ and $P\left\{ \left\Vert \mathbb{Z}\right\Vert _{\mathbb{T}%
}\leq M\left( X\right) \right\} \geq1/2$. We shall assume that $M\left(
X\right) $ is finite. \smallskip

\noindent\textbf{Borell's inequality.} \textit{For all $z>0$, 
\begin{equation}
P\left\{ \left\vert \left\Vert \mathbb{Z}\right\Vert _{\mathbb{T}}-E\left(
\left\Vert \mathbb{Z}\right\Vert _{\mathbb{T}}\right) \right\vert >z\right\}
\leq2\exp\left( -\frac{z^{2}}{2\sigma_{\mathbb{T}}^{2}\left( \mathbb{Z}%
\right) }\right) .  \label{bor}
\end{equation}
}

\subsubsection{Application of the Landau--Shepp Theorem}

We shall be using the
following version of the Landau and Shepp [LS] \cite{LandauShepp} theorem (also see 
Sat\^{o} \cite{Sato}, Theorem 2.5 of Marcus and Shepp \cite{MS} and Proposition A.2.3 in
\cite{VaartWellner}): \medskip

\noindent \textbf{Theorem [LS]} \textit{Let $X_{t}$, $t\in T,$ be a
real valued separable Gaussian process such that w.p.~$1$, $\sup_{t\in
T}\left\vert X_{t}\right\vert <\infty$, then for any $0<\beta <1/\left( 2\sigma 
^{2}\right) $, where $\sigma^{2}=\sup_{t\in T}\mathop{Var}\left( X_{t}\right)$, 
for all $y$ sufficiently large 
\begin{equation}
P\left\{ \sup_{t\in T}\left\vert X_{t}\right\vert >y\right\} <\exp \left(
-\beta y^{2}\right) .  \label{LST}
\end{equation}
}

Recall the definition of $L$ in (\ref{MC}). Since $L$ is finite, w.p.~$1$,
we can apply the Landau and Shepp theorem to infer that for appropriate
constants $C>0$ and $D>0$, for all $t>0,$%
\begin{equation}
P\left\{ L>t\right\} \leq C\exp \left( -Dt^{2}\right) .  \label{FLS1}
\end{equation}

\subsection{Four maximal inequalities}
\label{four}

For the following inequalities recall the mean zero Gaussian process $G$
with covariance function defined in (\ref{EG}). Inequalities 1 and 2 are
required for the proof of Proposition \ref{prop:2}, and Inequalities 1A and
2A are needed in the proofs of Theorems 1 and 2.\smallskip

\noindent \textbf{Inequality 1.} \textit{For all $0<\varrho <\infty$ and $\delta 
>0$ we have for some constant $\mu (\delta)$ and all $z>0$
\begin{equation}
P\left\{ \sup_{(t,x)\in \lbrack 0,\varrho ]\times \mathbb{R}}t^{\delta
}\left\vert G\left( t,x\right) \right\vert >\varrho ^{\delta }2^{\delta }\mu
\left( \delta \right) +z\right\} \leq 2\exp \left( -\frac{z^{2}\varrho
^{-2\delta }}{2^{2\delta +1}}\right)  \label{ginq}
\end{equation}%
and for each $n\geq 1$ and for $t^{\delta }G^{\left( 1\right)
}\left( t,x\right) ,\dots ,t^{\delta }G^{\left( n\right) }\left( t,x\right) $
i.i.d.~$t^{\delta }G\left( t,x\right)$
\begin{equation}
P\left\{ \max_{1\leq m\leq n}\sup_{(t,x)\in \lbrack 0,\varrho ]\times 
\mathbb{R}}\left\vert \frac{1}{\sqrt{n}}\sum_{i=1}^{m}t^{\delta }G^{\left(
i\right) }\left( t,x\right) \right\vert >\varrho ^{\delta }2^{\delta }\mu
\left( \delta \right) +z\right\} \leq 4\exp \left( -\frac{z^{2}\varrho
^{-2\delta }}{2^{2\delta +1}}\right) .  \label{ginq1}
\end{equation}%
}

\noindent \textit{Proof} Define for any integer $k\geq 0$, 
\begin{equation*}
\mathcal{T}_{k}=\left[ 2^{-k},2^{-k+1}\right] \times \mathbb{R}.
\end{equation*}%
Theorem 5 in \cite{KKZ} implies that, w.p.~1, for each integer $k$, 
\begin{equation}
\sup \left\{ \left\vert G\left( t,x\right) \right\vert :\left( t,x\right)
\in \mathcal{T}_{k}\right\} <\infty .  \label{finite}
\end{equation}%
Notice that for any $k\geq 0$ 
\begin{equation*}
\sup \left\{ \left\vert G\left( t,x\right) \right\vert :\left( t,x\right)
\in \mathcal{T}_{k}\right\} \overset{\mathrm{D}}{=}\sup \left\{ \left\vert
G\left( t,x\right) \right\vert :\left( t,x\right) \in \mathcal{T}%
_{0}\right\} .
\end{equation*}%
Furthermore (\ref{finite}) and separability of $G\left( t,x\right) $ permits
us to apply the Landau--Shepp theorem (see (\ref{LST})) to get 
\begin{equation*}
\mu _{0}:=E\left( \sup \left\{ \left\vert G\left( t,x\right) \right\vert
:\left( t,x\right) \in \mathcal{T}_{0}\right\} \right) <\infty .
\end{equation*}%
Thus for any integer $K$%
\begin{equation*}
E\left( \sup \left\{ t^{\delta }\left\vert G\left( t,x\right) \right\vert
:\left( t,x\right) \in \left[ 0,2^{-K}\right] \times \mathbb{R}\right\}
\right)
\end{equation*}%
\begin{equation*}
\leq \mu _{0}\sum_{k=K}^{\infty }2^{-\delta k}=2^{-\delta K}\mu _{0}/\left(
1-2^{-\delta }\right) =:2^{-\delta K}\mu \left( \delta \right) .
\end{equation*}%
This implies that, w.p.~$1$, 
\begin{equation*}
\sup \left\{ t^{\delta }\left\vert G\left( t,x\right) \right\vert :\left(
t,x\right) \in \left[ 0,2^{-K}\right] \times \mathbb{R}\right\} <\infty .
\end{equation*}%
Also 
\begin{equation*}
\sup \left\{ \mathop{Var}\left( t^{\delta }G\left( t,x\right) \right)
:\left( t,x\right) \in \left[ 0,2^{-K}\right] \times \mathbb{R}\right\} \leq
2^{-2\delta K}.
\end{equation*}%
Applying Borell's inequality (\ref{bor}) with $\mathbb{Z}\left( t,x\right)
=t^{\delta }G\left( t,x\right) $, $\mathbb{T}=\left[ 0,2^{-K}\right] \times 
\mathbb{R}$, $E\left( \left\Vert \mathbb{Z}\right\Vert _{\mathbb{T}}\right)
\leq 2^{-\delta K}\mu \left( \delta \right) $ and $\sigma _{\mathbb{T}%
}^{2}\left( \mathbb{Z}\right) \leq 2^{-2\delta K}$, we get for all $z>0$ and
integers $K$ 
\begin{equation*}
P\left\{ \sup_{(t,x)\in \lbrack 0,2^{-K}]\times \mathbb{R}}t^{\delta
}\left\vert G\left( t,x\right) \right\vert >2^{-\delta K}\mu \left( \delta
\right) +z\right\} \leq 2\exp \left( -\frac{z^{2}2^{2\delta K}}{2}\right) .
\end{equation*}%
Choose any $0<\varrho <\infty $ and integer $K$ such that $2^{-K}\geq
\varrho >2^{-K-1}.$ We see that 
\begin{equation*}
2^{K+1}>\varrho ^{-1}\geq 2^{K}\geq \varrho ^{-1}/2.
\end{equation*}%
Hence $\left[ 0,\varrho \right] \times \mathbb{R}\subset \left[ 0,2^{-K}%
\right] \times \mathbb{R}$. Therefore 
\begin{equation*}
\begin{split}
& P\left\{ \sup \left\{ t^{\delta }\left\vert G\left( t,x\right) \right\vert
:\left( t,x\right) \in \left[ 0,\varrho \right] \times \mathbb{R}\right\}
>\varrho ^{\delta }2^{\delta }\mu \left( \delta \right) +z\right\} \\
& \leq P\left\{ \sup \left\{ t^{\delta }\left\vert G\left( t,x\right)
\right\vert :\left( t,x\right) \in \left[ 0,2^{-K}\right] \times \mathbb{R}%
\right\} >2^{-\delta K}\mu \left( \delta \right) +z\right\} \\
& \leq 2\exp \left( -\frac{z^{2}2^{2\delta K}}{2}\right) \leq 2\exp \left( -%
\frac{z^{2}\varrho ^{-2\delta }}{2^{2\delta +1}}\right) .
\end{split}%
\end{equation*}%
Inequality (\ref{ginq1}) follows from L\'{e}vy's inequality (see Proposition
A.1.2 in van der Vaart and Wellner \cite{VaartWellner}) along with
separability of the Gaussian process $t^{\delta }G\left( t,x\right) .$\hfill 
$\square \medskip $

\noindent \textbf{Inequality 1A.} \textit{For all $0<\gamma <1<T<\infty$ we have 
for some constant $\mu$  and all $z>0$
\begin{equation}
P\left\{ \sup_{(t,x)\in \mathcal{T}\left( \gamma \right) }\left\vert G\left(
t,x\right) \right\vert >\mu +z\right\} \leq 2\exp \left( -\frac{z^{2}}{2}%
\right)  \label{ginqA}
\end{equation}%
and for each $n\geq 1$ and $G^{\left( 1\right) }\left(
t,x\right) ,\dots ,G^{\left( n\right) }\left( t,x\right)$ i.i.d.~$G\left( t,x\right)$
\begin{equation}
P\left\{ \max_{1\leq m\leq n}\sup_{((t,x)\in \mathcal{T}\left( \gamma
\right) }\left\vert \frac{1}{\sqrt{n}}\sum_{i=1}^{m}G^{\left( i\right)
}\left( t,x\right) \right\vert >\mu +z\right\} \leq 4\exp \left( -\frac{z^{2}%
}{2}\right) .  \label{ginqA1}
\end{equation}%
}

\noindent \textit{Proof} Theorem 5 in \cite{KKZ} implies that, w.p.~1, 
\begin{equation}
\sup \left\{ \left\vert G\left( t,x\right) \right\vert :\left( t,x\right)
\in \mathcal{T}\left( \gamma \right) \right\} <\infty .  \label{finite1}
\end{equation}%
Furthermore (\ref{finite1}) permits us to apply the Landau--Shepp theorem to
get 
\begin{equation*}
\mu :=E\left( \sup \left\{ \left\vert G\left( t,x\right) \right\vert :\left(
t,x\right) \in \mathcal{T}\left( \gamma \right) \right\} \right) <\infty .
\end{equation*}%
Also 
\begin{equation*}
\sup \left\{ \mathop{Var}\left( G\left( t,x\right) \right) :(t,x)\in 
\mathcal{T}\left( \gamma \right) \right\} \leq 1.
\end{equation*}%
Applying Borell's inequality (\ref{bor}) with $\mathbb{Z}\left( t,x\right)
=G\left( t,x\right)$, $\mathbb{T}=\mathcal{T}\left( \gamma \right)$, $%
E\left( \left\Vert \mathbb{Z}\right\Vert _{\mathbb{T}}\right) =\mu$ and $%
\sigma_{\mathbb{T}}^{2}\left( \mathbb{Z}\right) \leq 1$, we get for all $z>0$
\begin{equation*}
P\left\{ \sup_{(t,x)\in \mathcal{T}\left( \gamma \right) }\left\vert G\left(
t,x\right) \right\vert >\mu +z\right\} \leq 2\exp \left( -\frac{z^{2}}{2}
\right).
\end{equation*}%
Inequality (\ref{ginqA1}) follows from L\'{e}vy's inequality and
separability of the Gaussian process $G\left( t,x\right)$. \hspace*{1pt}%
\hfill $\square$ \medskip

In Inequalities 2 and 2A for $g\in \mathcal{C}\left[ 0,T\right]$,
\begin{equation*}
h_{t,x}\left( g\right) =1\left\{ g\left( t\right) \leq x,\text{ }g\in 
\mathcal{C}_{\infty }\right\},
\end{equation*}%
where $\mathcal{C}_{\infty }$ is defined as in (10).

\noindent \textbf{Inequality 2.} \textit{For all $0<\varrho < T/2$ and $%
\delta >0$ we have for some $E(\delta )$ and for suitable finite positive
constants $A, A_{1}>0$, for all $z>0$ 
\begin{equation}
\begin{split}
& P\left\{ \max_{1\leq m\leq n}\sup_{(t,x)\in \lbrack 0,\varrho ]\times 
\mathbb{R}}|\sqrt{m}t^{\delta }\alpha _{m}\left( h_{t,x}\right) |>\sqrt{n}%
A\left( E(\delta )2^{\delta }\varrho ^{\delta }+z\right) \right\} \\
& \leq 2\left\{ \exp \left( -z^{2}A_{1}\left( 2\varrho \right) ^{-2\delta
}\right) +\exp \left( -z\sqrt{n}A_{1}\left( 2\varrho \right) ^{-\delta
}\right) \right\} .
\end{split}
\label{H6}
\end{equation}
}

Note, in particular, Inequality 2 implies that for all $\lambda >1$ there is
a $d>1$ such that 
\begin{equation}
P\left\{ \sup \left\{ |t^{\delta }\alpha _{n}\left( h_{t,x}\right) |:\left(
t,x\right) \in \left[ 0,\varrho \right] \times \mathbb{R}\right\} \geq
d\varrho ^{\delta }\sqrt{\log n}\right\} <n^{-\lambda }.  \label{videlta}
\end{equation}

\noindent \textit{Proof} For any $k\geq 1$ and $g\in \mathcal{C}[0, T] $, 
let%
\begin{equation*}
g_{k}\left( t\right) =2^{kH}g\left( t2^{-k}\right) ,\ t\in [ 0,T] ,
\end{equation*}%
and for any $k\geq 1$, $t\in [ 0,T ] $, $x\in \mathbb{R}$
and $g\in \mathcal{C} [ 0, T]$ set
\begin{equation*}
h_{t,x,k}\left( g\right) =h_{t,x}\left( g_{k}\right) =1\{ g_{k}(t) \leq x, g_{k}\in 
\mathcal{C}_{\infty } \} .
\end{equation*}%
Clearly w.p.~$1$%
\begin{equation*}
\sup_{(t,x)\in \mathcal{T}_{k}}\left\vert
\sum_{i=1}^{n}\epsilon_{i}h_{t,x}(B_{i}) \right\vert =\sup_{(t,x)\in 
\mathcal{T}_{0}}\left\vert \sum_{i=1}^{n}\epsilon
_{i}h_{t,x,k}(B_{i})\right\vert.
\end{equation*}%
Moreover, since 
\begin{equation*}
\left\{ B_{j}\right\} _{j\geq 1}\overset{\mathrm{D}}{=}\left\{
2^{kH}B_{j}\left( \cdot /2^{k}\right) \right\} _{j\geq 1}\text{,}
\end{equation*}%
we see that%
\begin{equation*}
\sup_{(t,x)\in \mathcal{T}_{k}}\left\vert \sum_{i=1}^{n}\epsilon
_{i}h_{t,x,k}(B_{i})\right\vert \overset{\mathrm{D}}{=}\sup_{(t,x)\in 
\mathcal{T}_{0}}\left\vert \sum_{i=1}^{n}\epsilon
_{i}h_{t,x}(B_{i})\right\vert
\end{equation*}%
and thus 
\begin{equation}
E\sup_{(t,x)\in \mathcal{T}_{k}}\left\vert \sum_{i=1}^{n}\epsilon
_{i}h_{t,x}(B_{i})\right\vert =E\sup_{(t,x)\in \mathcal{T}_{0}}\left\vert
\sum_{i=1}^{n}\epsilon _{i}h_{t,x}(B_{i})\right\vert .  \label{ID}
\end{equation}%
We readily see by inequality (\ref{H5A}) 
\begin{equation*}
E\sup_{(t,x)\in \mathcal{T}_{0}}\left\vert \sum_{i=1}^{n}\epsilon
_{i}h_{t,x}(B_{i})\right\vert \leq 2\sqrt{n}E\left\Vert v_{n}\right\Vert _{%
\mathcal{T}_{0}}+ \sqrt{n},
\end{equation*}%
which by (\ref{M}) is $\leq 2\left( M\left( 1,2,H\right) +1\right) \sqrt{n}%
=:E_{0}\sqrt{n}.$ Thus 
\begin{equation}
E\sup_{(t,x)\in \mathcal{T}_{0}}\left\vert \sum_{i=1}^{n}\epsilon
_{i}h_{t,x}(B_{i})\right\vert \leq E_{0}\sqrt{n}.  \label{33}
\end{equation}%
Next, for all $\delta >0$ with $K$ an integer such that $2^{-K}\geq \varrho >2^{-K-1}$
\begin{equation*}
E\sup \left\{ \left\vert t^{\delta }\sum_{i=1}^{n}\epsilon
_{i}h_{t,x}(B_{i})\right\vert :0\leq t\leq 2^{-K},x\in \mathbb{R}\right\}
\end{equation*}%
is by (\ref{ID}) and (\ref{33}) 
\begin{equation}
\leq \sum_{k=K}^{\infty }2^{-k\delta }E\sup_{(t,x)\in \mathcal{T}%
_{k}}\left\vert \sum_{i=1}^{n}\epsilon _{i}h_{t,x}(B_{i})\right\vert \leq
E\left( \delta \right) 2^{-K\delta }\sqrt{n\,},  \label{E3}
\end{equation}%
where $E\left( \delta \right) =E_{0}/\left( 1-2^{-\delta }\right) .$
\smallskip

Let%
\begin{equation*}
\mathcal{H}\left( \delta ,K\right) =\left\{ t^{\delta }h_{t,x}:\left(
t,x\right) \in \left[ 0,2^{-K}\right] \times \mathbb{R}\right\} .
\end{equation*}%
From (\ref{E3}) we get 
\begin{equation}
E\sup \left\{ \left\vert \sum_{i=1}^{n}\epsilon _{i}g(B_{i})\right\vert
:g\in \mathcal{H}\left( \delta ,K\right) \right\} \leq E\left( \delta
\right) 2^{-K\delta }\sqrt{n\,}.  \label{E2}
\end{equation}%
Also observe that each $g\in \mathcal{H}\left( \delta ,K\right)$ satisfies $%
|g|\leq 2^{-K\delta }$. Applying Talagrand's inequality (\ref{tal}) with $%
M=2^{-K\delta }$, $\sigma _{\mathcal{H}\left( \delta ,K\right) }^{2} =
2^{-2K\delta}$ and the bound (\ref{E2}), we get that for any $\delta >0$ we
have for suitable finite positive constants $A,A_{1}>0$, for all $z>0$ 
\begin{equation}
\begin{split}
& P\left\{ \max_{1\leq m\leq n}||\sqrt{m}\alpha _{m}||_{\mathcal{H}\left(
\delta ,K\right) }\geq \sqrt{n}A(E\left( \delta \right) 2^{-K\delta
}+z)\right\} \\
& \leq 2(\exp (-z^{2}A_{1}2^{2K\delta })+\exp (-z\sqrt{n}A_{1}2^{K\delta })).
\end{split}
\label{H7}
\end{equation}

Inequality (\ref{H6}) follows from inequality (\ref{H7}). To see this choose
any $0<\varrho <T/2 $ and integer $K$ such that $2^{-K}\geq \varrho
>2^{-K-1}.$ We see that 
\begin{equation*}
2^{K+1}>\varrho ^{-1}\geq 2^{K}\geq \varrho ^{-1}/2.
\end{equation*}%
Hence $\left\{ t^{\delta }h_{t,x}:\left( t,x\right) \in \left[ 0,\varrho %
\right] \times \mathbb{R}\right\} \subset \mathcal{H}\left( \delta ,K\right)$, and 
\begin{equation*}
\begin{split}
& P\left\{ \max_{1\leq m\leq n}\sup_{(t,x)\in \lbrack 0,\varrho ]\times 
\mathbb{R}}|\sqrt{m}t^{\delta }\alpha _{m}\left( h_{t,x}\right) |\geq \sqrt{n%
}A(E\left( \delta \right) 2^{\delta }\varrho ^{\delta }+z)\!\right\} \\
& \leq P\left\{ \max_{1\leq m\leq n}||\sqrt{m}\alpha _{m}||_{\mathcal{H}%
\left( \delta ,K\right) }\geq \sqrt{n}A(E\left( \delta \right) 2^{-K\delta
}+z)\right\} \\
& \leq 2(\exp (-z^{2}A_{1}2^{2K\delta })+\exp (-z\sqrt{n}A_{1}2^{K\delta }))
\\
& \leq 2\left\{ \exp \left( -z^{2}A_{1}(2\varrho )^{-2\delta }\right) +\exp
\left( -z\sqrt{n}A_{1}(2\varrho )^{-\delta }\right) \right\} .
\end{split}%
\end{equation*}%
\qed \pagebreak

\noindent \textbf{Inequality 2A.} \textit{For all $0<\gamma <1<T<\infty $,
we have for some some $L(\gamma ,T)$ and all $z>0$ for suitable finite
positive constants $A, A_{1}>0$, for all $z>0$ 
\begin{equation}
\begin{split}
& P\left\{ \max_{1\leq m\leq n}\sup_{(t,x)\in \mathcal{T}\left( \gamma
\right) }|\sqrt{m}\alpha _{m}\left( h_{t,x}\right) |\geq \sqrt{n}A\left(
L(\gamma ,T)+z\right) \right\} \\
& \leq 2\left\{ \exp \left( -z^{2}A_{1}\right) +\exp \left( -z\sqrt{n}%
A_{1}\right) \right\} .
\end{split}
\label{in2A}
\end{equation}
} 

\noindent \textit{Proof} We see by inequality (\ref{H5A}) 
\begin{equation*}
E\sup_{(t,x)\in \mathcal{T}\left( \gamma \right) }\left\vert
\sum_{i=1}^{n}\epsilon _{i}h_{t,x}(B_{i})\right\vert \leq 2\sqrt{n}%
E\left\Vert v_{n}\right\Vert _{\mathcal{T}\left( \gamma \right) }+ \sqrt{n},
\end{equation*}%
which by (\ref{M}) is $\leq 2\left( M\left( \gamma ,T,H\right) +1\right) 
\sqrt{n}=:L(\gamma ,T)\sqrt{n}.$ Thus 
\begin{equation}
E\sup_{(t,x)\in \mathcal{T}\left( \gamma \right) }\left\vert
\sum_{i=1}^{n}\epsilon _{i}h_{t,x}(B_{i})\right\vert \leq L(\gamma ,T)\sqrt{n%
}.  \label{Lg}
\end{equation}%
Applying Talagrand's inequality (\ref{tal}) with $M=1$, $\sigma _{\mathcal{F}%
_{\left( \gamma ,T\right) }}^{2}=1$ and the bound (\ref{Lg}), give (\ref%
{in2A}).\newline
\hspace*{1pt}\hfill $\square $\medskip

\begin{remark}
\label{rem:7} Actually, to apply Talagrand's inequality in the proofs of
Inequalities 2 and 2A, as it is stated in (\ref{tal}), the classes of
functions $\mathcal{H}\left( \delta ,K\right) $ and $\mathcal{F}_{\left(
\gamma ,T\right) }$ should be pointwise measurable. Here we shall discuss
how to take care of this detail in the proof of Inequality 2. A similar
discussion works for the proof of Inequality 2A.

For any $k\geq 1$ let 
\begin{equation*}
\mathcal{H}\left( \delta ,K,k\right) =\left\{ g1\left\{ g\in \mathcal{C}%
\left( k\right) \right\} :g\in \mathcal{H}\left( \delta ,K\right) \right\},
\end{equation*}
where $\mathcal{C}\left( k\right) $ is defined as in (\ref{cck}).
The class $\mathcal{H}\left( \delta ,K,k\right) $ is pointwise measurable.
Applying Talagrand's inequality we get with $M=2^{-K\delta }$ and $\sigma _{%
\mathcal{H}\left( \delta ,K,k\right) }^{2}=2^{-2K\delta }$ 
\begin{equation*}
\begin{split}
& P\left\{ \max_{1\leq m\leq n}||\sqrt{m}\alpha _{m}||_{\mathcal{H}\left(
\delta ,K,k\right) }\geq A\left( E\bigg\Vert \sum_{i=1}^{n}\epsilon
_{i}g(B_{i})\bigg\Vert _{\mathcal{H}\left( \delta ,K,k\right) }+t\right)
\right\} \\
& \leq 2\exp \left( -\frac{2^{2K\delta }A_{1}t^{2}}{n}\right) +2\exp \left(
-2^{K\delta }A_{1}t\right) .
\end{split}%
\end{equation*}%
Obviously by the Wang \cite{Wang} result (\ref{W}), w.p.~$1$, $B \in \cup
_{k=1}^{\infty }\mathcal{C}\left( k\right)$. Therefore, w.p.~$1$, for any $%
n\geq 1$, $B_{1},\dots ,B_{n}$, i.i.d.~$B$ there exists a $k\geq 1$ such
that uniformly in $\left( t,x\right) \in \left[ 0,\varrho \right] \times 
\mathbb{R}$, $h_{t,x}^{\left( k\right) }\left( B_{i} \right) =h_{t,x}\left(
B_{i} \right) $ and $t^{\delta}h_{t,x}^{\left( k\right) } \left( B_{i}
\right) =t^{\delta }h_{t,x}\left( B_{i}\right)$, for $i=1,\dots,n$. This
says that, w.p.~$1$, for any $n\geq 1$, there exists a $k\geq 1,$ such that
uniformly in $\left( t,x\right) \in \left[ 0,\varrho \right] \times \mathbb{R%
}$ 
\begin{equation*}
\frac{1}{\sqrt{n}}\sum_{i=1}^{n}t^{\delta }h_{t,x}\left( B_{i}\right)
1\left\{ B_{i}\notin \mathcal{C}\left( k\right) \right\} =0.
\end{equation*}%
Furthermore%
\begin{equation*}
\sup_{\left( t,x\right) \in \left[ 0,\varrho \right] \times \mathbb{R}}\frac{%
1}{\sqrt{n}}\sum_{i=1}^{n}t^{\delta} Eh_{t,x}\left( B_{i} \right) 1\left\{
B_{i}\notin \mathcal{C}\left( k\right) \right\} \leq \sqrt{n}\varrho
^{\delta } P \left\{ B\notin \mathcal{C}\left( k\right) \right\},
\end{equation*}%
which converges to zero for each fixed $n\geq 1$, as $k\rightarrow \infty $.
By passing to the limit, as $k\rightarrow \infty $, we get for any $\delta
>0 $ and $t>0$%
\begin{equation*}
\begin{split}
& P\left\{ \max_{1\leq m\leq n}||\sqrt{m}\alpha _{m}||_{\mathcal{H}\left(
\delta ,K\right) }\geq A\left( E\left\Vert \sum_{i=1}^{n}\epsilon
_{i}g(B_{i})\right\Vert _{\mathcal{H}\left( \delta ,K\right) }+t\right)
\right\} \\
& \leq 2\exp \left( -\frac{2^{2K\delta }A_{1}t^{2}}{n}\right) +2\exp \left(
-2^{K\delta }A_{1}t\right) .
\end{split}%
\end{equation*}%
Similarly one can argue the validity of the Talagrand inequality using the
index class $\mathcal{F}_{\left( \gamma ,T\right) }$.
\end{remark}

\noindent\textbf{Acknowledgement} The authors thank the Associate Editor for
a comment that led to Remark \ref{rem:c}. PK was partially supported by the
Hungarian Scientific Research Fund OTKA PD106181, by the European Union and
co-funded by the European Social Fund under the project
`Telemedicine-focused research activities on the field of Mathematics,
Informatics and Medical sciences' of project number 
T\'{A}MOP-4.2.2.A-11/1/KONV-2012-0073 and by a postdoctoral fellowship of the Alexander 
von Humboldt Foundation.

\end{document}